\documentclass[12pt]{article}

\usepackage[english]{babel}

\usepackage{amsmath,amssymb,amsthm}
\usepackage[mathscr]{eucal}
\usepackage{bbm}

\newcommand{\C}{\mathbb{C}}
\newcommand{\D}{\mathbb{D}}
\newcommand{\Disk}{\mathbb{D}}

\newcommand{\N}{\mathbb{N}}
\newcommand{\R}{\mathbb{R}}

\newcommand{\Z}{\mathbb{Z}}

\newcommand{\A}{\mathcal{A}}

\newcommand{\Dr}[1]{\mathcal{D}^{#1}}
\newcommand{\scd}{\mathcal{D}}
\newcommand{\F}{\mathcal{F}}

\newcommand{\Square}{\mathcal{S}}

\newcommand{\e}{\varepsilon}
\newcommand{\eps}{\varepsilon}
\renewcommand{\epsilon}{\varepsilon}

\renewcommand{\phi}{\varphi}

\newcommand{\LE}{\mathcal{L}}

\renewcommand{\d}{\mathrm{d}}
\newcommand{\n}{\mathbf{n}}

\renewcommand{\P}{\mathbb{P}} 

\newcommand{\E}{\mathbb{E}}

\newcommand{\bd}{\partial}  

\newcommand{\ds}{\displaystyle}

\newcommand{\rad}{\operatorname{rad}}
\newcommand{\inrad}{\operatorname{inrad}}
\newcommand{\dist}{\operatorname{dist}}
\newcommand{\diam}{\operatorname{diam}}
\newcommand{\sgn}{\operatorname{sgn}}

\DeclareMathAlphabet{\mathpzc}{OT1}{pzc}{m}{it}
\renewcommand{\Re}{\mathpzc{Re}}
\renewcommand{\Im}{\mathpzc{Im}}

\newcommand{\edgeset}[1]{\mathscr{E}_{#1}}
\newcommand{\univalent}{\mathscr{S}}

\newcommand{\defines}[1]{\textit{#1}\index{#1}}

\newcommand{\excdetfrac}{\Lambda}
\newcommand{\FominEvent}{\mathcal{C}}

\newcommand{\excur}{\mathcal{E}}
\newcommand {\p} {\partial}
\newcommand {\Half} {\mathbb{H}}
\newcommand{\rect} {\mathcal{R}}

\newtheorem{theorem}{Theorem}[section]
\newtheorem*{theorem*}{Theorem}
\newtheorem{proposition}[theorem]{Proposition}
\newtheorem{lemma}[theorem]{Lemma}

\newtheorem{corollary}[theorem]{Corollary}
\newtheorem*{corollary*}{Corollary}

\newtheorem*{conjecture*}{Conjecture}

\theoremstyle{remark}

\newtheorem*{remark}{Remark}

\theoremstyle{definition}
\newtheorem{definition}[theorem]{Definition}

\newtheorem{example}[theorem]{Example}

\author{Michael J.~Kozdron\footnote{University of Regina} \and
Gregory F.~Lawler\footnote {Cornell University (Research supported by the 
National Science Foundation.)}}

\title{Estimates of random walk exit probabilities
and application to loop-erased random walk}

\date{January 12, 2005} 

\begin{document}

\maketitle


\begin{abstract}
We prove an estimate for the probability that a simple random walk in a
simply connected subset $A \subset \Z^2$ starting on the boundary exits
$A$ at another specified boundary point.  The estimates are uniform
over all domains of a given inradius.  We apply these estimates to
  prove a conjecture  of S.~Fomin~\cite{Fom1} in 2001 concerning a relationship 
between crossing
probabilities of loop-erased random walk and Brownian motion.    
\end{abstract}



\section{Introduction}\label{Sect1}

In the last few years a number of results have been proved about 
scaling limits of two-dimensional lattice systems in statistical
mechanics.  Site percolation on the triangular lattice~\cite{Smirnov},
loop-erased random walk~\cite{LawSW9,Zhan}, uniform spanning trees~\cite{LawSW9}, and the harmonic explorer~\cite{SS} have all been
shown to have limits that can be described using the Schramm-Loewner
evolution.  In the last three cases, the proofs use a 
version of the well-known fact that simple random walk  
has a scaling limit of Brownian motion, which is
conformally invariant in two dimensions.  What is needed
is a strong version of this result which holds uniformly over
a wide class of domains where the errors do not depend on
the smoothness of the boundary.  In this paper, we present another
result of this type.  It differs from the lemmas in~\cite{LawSW9}
in two ways: we give explicit error bounds that show that
the error decays as a power of the ``inradius'' of the domain, and
the continuous domain that we compare to the discrete domain
is slightly different.

We give an application of our result to loop-erased walk by
proving a conjecture of S.~Fomin~\cite{Fom1} and giving a quick
derivation of a crossing exponent first proved by R.~Kenyon~\cite{Kenyon}.
Fomin showed that a certain crossing probability for loop-erased walk
can be given in terms of a determinant of hitting probabilities
for simple random walk.  Our estimate shows that this determinant
approaches a corresponding determinant for Brownian motion.  We then
estimate the determinant for Brownian motion to derive
the  crossing exponent.

We will start by discussing the main results, leaving some of the
precise definitions until Section~\ref{Sect2}.
The only Euclidean dimension that will concern us is $d=2$; consequently, 
we associate $\C \cong \R^2$ in the natural way. Points in the complex plane 
will be denoted by any of $w$, $x$, $y$, or $z$.  
A \defines{domain} $D \subset \C$ is an open and connected set.

Throughout this paper, $B_t$, $t\ge 0$, will denote a standard complex 
Brownian motion, and $S_n$, $n=0, 1, \ldots$, will 
denote two-dimensional simple random walk, 
both started at the origin unless otherwise noted. 
We write $B[0,t] := \{ z \in \C : B_s = z \, \text{ for some } \, 0 \le s \le t \}$,
 and $S[0,n] := [S_0, S_1, \ldots, S_n]$ 
for the set of lattice points visited by the random walk. 
We will generally use $T$ for stopping times for Brownian motion, 
and $\tau$ for stopping times for random walk. 
We write $\E^x$ and $\P^x$ for expectations and probabilities, 
respectively, assuming $B_0=x$ or $S_0=x$, as appropriate.

\subsection{Main results}
Let $\A^n$  denote the collection of simply connected subsets $A$ of $\Z^2$ such that $n \leq \inrad(A) \leq 2n$, i.e., such that
\[   n \leq \sup\{|z|: z \in \Z^2 \setminus A \} \leq 2n . \]
Associated to $A$ is a simply connected domain $\tilde A \subset \C$ which is obtained by identifying each lattice point in $A$ with the square of side one centred at that point. By the Riemann mapping theorem, there is a unique conformal transformation $f_A$ of $\tilde A$ onto the unit disk with $f_A(0) = 0, f_A'(0) > 0$. We let $\theta_A(x) := \arg(f_A(x))$.  We can extend $\theta_A$ to $\bd A$ in a natrual way.

If $x \in A$, $y \in \bd A$, let $h_A(x,y)$
be the probability that a simple random walk starting
at $x$ leaves $A$ at $y$.
If $x$, $y \in \bd A$,  let $h_{\bd A}(x,y)$ be the probability
that a simple random walk starting at $x$ takes its first step
into $A$ and then leaves $A$ at $y$.

\begin{theorem} \label{rwestimate}
 If $A \in \A^n$, then
\[   h_{\bd A}(x,y) = \frac{(\pi/2) \, h_A(0,x) \,
    h_A(0,y) }  { 1 - \cos(\theta_A(x) - \theta_A(y))}
   \; \left[1 + O\left(\frac{\log n}{n^{1/16} \, |\theta_A(x) -
  \theta_A(y)|}\right)\, \right] , \]
provided that $|  \theta_A(x) -
  \theta_A(y)| \geq n^{-1/16} \, \log^{2}n$.
\end{theorem}

In this theorem, and throughout this paper, we will use $O(\cdot)$
for uniform error terms that depend only on $n$.  For example, the
statement above is shorthand for the following: there is a constant
$c < \infty$ such that for all $A \in \A^n$ and all
$x,y \in \bd A$ with  $|\theta_A(x) -
  \theta_A(y)| \geq n^{-1/16} \, \log^{2}n$,
\[   \left| \, h_{\bd A}(x,y) - \frac{(\pi/2) \, h_A(0,x) \,
    h_A(0,y) }  { 1 - \cos(\theta_A(x) - \theta_A(y))}
\,  \right | \leq c \, \frac{(\pi/2) \, h_A(0,x) \,
    h_A(0,y) }  { 1 - \cos(\theta_A(x) - \theta_A(y))}
  \; \frac{\log n}{n^{1/16} \, |\theta_A(x) -
  \theta_A(y)|}. \]
When discussing $k$-fold
determinants, we will have error terms
that depend also on a positive integer $k$; we will write
these as $O_k(\cdot)$.

We do not believe
the error term   $O(n^{-1/16} \,\log n)$ is    optimal,
and we probably could have improved it slightly in this paper.
However, our methods are not strong enough to give the optimal
error term.  The importance of this result is that the error
is bounded uniformly over all simply connected domains and that
the error is in terms of a power of $n$. For domains with
``smooth'' boundaries, one can definitely improve the power
of $n$.

To help understand this estimate, one should
 consider $h_{\p A} (x,y)$
as having a ``local'' and a ``global'' part.  The local part,
which is very dependent on the structure of $A$ near $x$
and $y$, is represented by the $h_A(0,x) \, h_A(0,y)$ term.
The global part, which is
$[ 1 - \cos(\theta_A(x) - \theta_A(y))]^{-1}$,
 is the conformal invariant and  depends
only on the image of the points under the conformal
transformation of $\tilde{A}$ onto  the unit disk.
In contrast to the discrete case, Example~\ref{ex3} 
shows that the Brownian version of this result is exact.

As part of the proof, we also derive a uniform estimate for
$G_A(x)$, the expected number of visits to $x$ before leaving
$A$ of a simple random walk starting at $0$.  Let $a$
denote the potential kernel for two-dimensional simple
random walk.  It is known that there is a $k_0$ such
that 
\[   a(x) = \frac 2 \pi\, \log|x| + k_0 + O(|x|^{-2}) , \;\;\;\;
   x \rightarrow \infty . \]

\begin{theorem}  \label{greentheoremB}
 If $A \in \A^n$, then
\[ G_A(0) = -\frac 2 \pi \, f_A'(0) + k_0 + O(n^{-1/3}
 \, \log n). \]
Furthermore, if $x \neq 0$, then
\begin{equation} \label{mar29.eq1}  
G_A(x) = \frac{2}{\pi} \, g_A(x) + k_x + O(n^{-1/3} \,
  \log n)  . \end{equation}
where $g_A(x) := g_A(0,x) = -\log |f_A(x)| $ is the Green's function for Brownian
motion in $\tilde A$  and
\[   k_x :=  k_0 + \frac 2 \pi \, \log |x| - a(x) . \]
\end{theorem}

\subsection{Fomin's identity for loop-erased walk}

We briefly review 
the definition of the loop-erased random walk; see~\cite[Chapter 7]{LawlerGreen} and~\cite{LawKesten}
 for more
details. Since simple random walk in $\Z^2$ 
is recurrent, it is not possible to construct loop-erased random walk by 
erasing loops from an infinite walk.  However, the following loop-erasing 
procedure makes perfect sense since it assigns to each finite simple random walk
 path a self-avoiding walk.  Let $S := S[0,m] := [S_0, S_1, \ldots, S_m]$ be a simple
 random walk path of length $m$. We construct $\LE(S)$, the loop-erased
 part of $S$, recursively as follows. If $S$ is already self-avoiding, 
set $\LE(S)=S$.  Otherwise, let $s_0 = \max\{j : S_j=S_0\}$, and for $i > 0$,
 let $s_i = \max\{j : S_j = S_{s_{j-1}+1} \}$. If we let $n = \min\{i : s_i=m\}$, then 
$\LE(S) = [S_{s_0}, S_{s_1}, \ldots, S_{s_n}]$.

Suppose that $A \in \A^n$ and $x^1,\ldots,x^k,y^k,
\ldots,y^1$ are distinct points in $\bd A$, ordered 
counterclockwise. For $i=1, \ldots, k$, 
let $\LE^i = \LE(S^i)$ be the loop erasure of the 
path $[S^i_0 = x^i, S^i_1, \ldots, S^i_{\tau^i_A}]$, 
and let $\FominEvent = \FominEvent(x^1, \ldots, x^k, y^k, 
\ldots, y^1; A)$ be the event that both
\begin{equation}  \label{dec14.1}
 S^i_{\tau^i_A} = y^i, \quad i=1, \ldots, k,
\end{equation}
 and
\begin{equation}  \label{dec14.2}
 S^i[0, \tau^i_A]\cap (\LE^1 \cup \cdots \cup \LE^{i-1})=
\emptyset, \quad i=2,\ldots, k.
\end{equation}

The following theorem was proved in~\cite{Fom1} which relates a determinant of simple random walk probabilities to a ``crossing probability'' for loop-erased random walk.  

\begin{theorem}[Fomin]\label{fomintheorem}
If $\FominEvent$ is the event defined above, and 
\begin{equation*}
\mathbf{h}_{\bd A}(\mathbf{x},\mathbf{y}) := 
\begin{bmatrix}
h_{\bd A}(x^1, y^1) &\cdots & h_{\bd A}(x^1, y^k) \\
\vdots      &\ddots &\vdots       \\
h_{\bd A}(x^k, y^1) &\cdots & h_{\bd A}(x^k, y^k)
\end{bmatrix},
\end{equation*}
where $\mathbf{x} = (x^1, \ldots, x^k)$, $\mathbf{y}  =(y^1, \ldots, y^k)$, 
then $\P(\FominEvent) = \det \mathbf{h}_{\bd A}  (\mathbf{x},\mathbf{y})$.
\end{theorem}

This is a special case of an identity that Fomin established for
 general discrete stationary Markov processes.  In his paper,
he  made the following conjecture.
\begin{quote}
In order for the statement of 
Theorem~\ref{fomintheorem} to make sense, 
the Markov process under consideration does not have to be 
discrete\dots. The proofs can be obtained by passing to a limit 
in the discrete approximation.  The same limiting procedure can be 
used to justify the well-definedness of the quantities involved; 
notice that in order to define a continuous analogue of Theorem~\ref{fomintheorem}, 
we do not need the notion of loop-erased Brownian motion.  
Instead, we discretize the model, compute the probability, 
and then pass to the limit.  One can further extend these results 
to densities of the corresponding hitting distributions. Technical details are omitted.
\end{quote}

With Theorem~\ref{rwestimate}, we have taken care of the
``technical details'' in the case of simply
connected planar domains.  Note that
\begin{equation}\label{conditionalFomin}
   \det \left[\frac{h_{\bd A}(x^j,y^l)}{h_{\bd A}(x^j,y^j)}\right]_{1 \leq 
   j,l \leq k} 
=\frac{ \ds \det\mathbf{h}_{\bd A}  (\mathbf{x},\mathbf{y}) }{\ds \prod_{j=1}^k h_{\bd A}(x^j,y^j)} 
\end{equation}
represents the conditional probability that~(\ref{dec14.2})
holds given~(\ref{dec14.1}) holds.  Suppose $D$
is a smooth
Jordan domain, and  that $x^1,\ldots,x^k,y^k,\ldots,y^1$
are distinct points on $\bd D$ ordered counterclockwise.
The ``Brownian motion'' analogue of the determinant~(\ref{conditionalFomin}) is
\begin{equation}\label{conditionalFominBM}
 \excdetfrac_{D}(x^1,\ldots,x^k,y^k,\ldots,y^1)
              :=  \det \left[\frac{H_{\bd D}(x^j,y^l)}
              {H_{\bd D}(x^j,y^j)}\right]_{1 \leq 
   j,l \leq k} 
= \frac{\det \mathbf{H}_{\bd D} (\mathbf{x},\mathbf{y}) }
  {\ds \prod_{j=1}^k H_{\bd D}(x^j,y^j)},  
\end{equation}
where $H_{\bd D}(z,w)$ denotes the excursion Poisson kernel.
In the case $D = \Disk$, if $z = e^{i \theta},
  w = e^{i \theta'}$, then
\[   H_{\bd \D}(z,w) = \frac 1 \pi \, \frac 1{|w-z|^2}
   = \frac 1 {2 \pi} \frac{1}{1- \cos(\theta' - \theta)}. \]
Conformal covariance of the excursion Poisson kernel
shows that
\begin{align*}
  \excdetfrac_{ D}(x^1,\ldots,x^k,y^k,\ldots,y^1)
              & =  \excdetfrac_{\Disk}(f(x^1),\ldots,f(x^k),
f(y^k),\ldots,f(y^1)) \\
&=  \det\left[\frac{1- \cos(\theta_D(x^j)
   - \theta_D(y^j))}{1 - \cos(\theta_D(x^j) -
  \theta_D(y^l))}\right]_{1 \leq j,l \leq k}, 
\end{align*}
where $f$ is a conformal transformation of $D$ onto
$\Disk$ and $\theta_D(z) := \arg (f(z))$. 

\begin{corollary}  Suppose 
$A \in \A^n$ and $x^1,\ldots,x^k,y^k,
\ldots,y^1$ are distinct points in $\bd A$ ordered
counterclockwise. Let
\[ m  = \min\{  \; 
  |\theta_A(x^1) - \theta_A(y^1)| , \;  |\theta_A(x^k) -
   \theta_A(y^k)| \; \}. \]
If $m \geq n^{-1/16} \, \log^{2}n$, then
\begin{equation}  \label{dec13.1}
  \det \left[\frac{h_{\bd A}(x^j,y^l)}{h_{\bd A}(x^j,y^j)}\right]_{1 \leq 
   j,l \leq k}
=    \det\left[\frac{1- \cos(\theta_A(x^j)
   - \theta_A(y^j))}{1 - \cos(\theta_A(x^j) -
  \theta_A(y^l))}\right]_{1 \leq j,l \leq k} +
   O_k\left(\frac{\log n}{n^{1/16} \,m^{2k+1}  
  } \right).
       \end{equation}              
\end{corollary}

\begin{proof}  Theorem~\ref{rwestimate} gives
 \[
\det \left[\frac{h_{\bd A}(x^j,y^l)}{h_{\bd A}(x^j,y^j)}\right] 
=   \det\left[\frac{1- \cos(\theta_A(x^j)
   - \theta_A(y^j))}{1 - \cos(\theta_A(x^j) -
  \theta_A(y^l))} \, [1 + O\left(\frac{\log n}{m \, n^{1/16}}\right)]\right].
\]
But, if $|\delta_{j,l}| \leq \epsilon$, multilinearity of the determinant
and the estimate $\det[b_{j,l}]  \leq k^{k/2}
  \, [\sup |b_{j,l}|]^k$
shows that
\begin{equation*}
   \left| \, \det[b_{j,l}(1 + \delta_{j,l})]
    - \det[b_{j,l}] \, \right| \leq [(1 + \epsilon)^k - 1]
               k^{k/2} \,  [\sup|b_{j,l}|]^k.  \qedhere
\end{equation*}
\end{proof}

Using the corollary, we know that we can approximate   the determinant
for random walks, and hence the probability of the crossing event $\FominEvent$,
in terms of the corresponding quantity for Brownian motion, at least
for simply connected domains.  We will consider
the asymptotics of   $\excdetfrac_{D}(x^1,\ldots,x^k,y^k,\ldots,y^1)$
   when $x^1,\ldots,x^k$ get close and
$y^1,\ldots,y^k$ get close.
Since this quantity is a conformal invariant, we may assume that
$D = \rect_L$, where
\[       \rect_L = \{\,z: \, 0 < \Re(z) < L, \,  0 < \Im(z) < \pi\,\}, \]
and $x^j = iq_j, y^j = L + i q_j'$, where $0 < q_k < \cdots < q_1 < \pi$
and $0 < q_k' < \cdots < q_1' < \pi$.   

\begin{proposition}  \label{propdec13.1}
As $L \rightarrow
\infty$,
\begin{align*} 
 \excdetfrac_{\rect_L}(iq_1,\ldots, &iq_k, L + iq_k',\ldots,
   L + i q_1')\\ 
  &= k! \, \frac{\det[\sin(lq_j)]_{1\le j,l \le k}  \;
          \det[\sin(lq_j')]_{1\le j,l \le k}  }
     {\ds \prod_{j=1}^k \sin(q_j) \, \sin(q_j')}
\,  
     e^{-k(k-1)L/2} + O_k(e^{-k(k+1)L/2}) .
\end{align*}
\end{proposition}

This crossing exponent $k(k-1)/2$ was first proved by Kenyon~\cite{Kenyon}
for loop-erased walk.  We describe this result in the framework of
Brownian excursion measure.

\subsection{Outline of the paper}  In the first five subsections
of Section~\ref{Sect2}, we review facts about random walk, Brownian motion,
and conformal mapping that we will need.   In the remaining subsections, we review 
Brownian excursion measure, define the analogue of the Fomin determinant
for excursion measure, and then derive Proposition~\ref{propdec13.1}.
In Section~\ref{Sect3}, we begin with a brief review of strong approximation, before proving Theorem~\ref{greentheoremB} in the second subsection. The final two subsections contain the proof of the other main result, Theorem~\ref{rwestimate}.


\section{Background, notation, and preliminary results}\label{Sect2}

\noindent In this section we review some basic material 
that will be needed in subsequent parts and standardize
 our notation.  Almost all of the complex analysis is well-known, 
and may be found in a variety of sources; we prove several elementary 
results, but often refer the reader to the literature for details. 
The material on the \emph{excursion Poisson kernel} is not difficult, 
but these results are not widespread.

\subsection{Simply connected subsets of $\C$ and $\Z^2$}\label{scs}

We will use $D$ to denote {\em domains}, i.e.,
 open connected subsets of $\C$.
We write $\D := \{z \in \C : |z|< 1\}$ to denote the open unit disk, 
and $\Half := \{z \in \C: \Im(z) > 0\}$ to denote the
upper half plane. 
 An analytic, univalent (i.e, one-to-one)
 function\footnote{For an analytic function $f$, 
$f'(z_0) \neq 0$ if and only if $f$ is locally univalent at $z_0$. 
 However, we will not be concerned with local univalence.} 
 is called a \defines{conformal mapping}.
 We say that $f:D \to D'$ is a \defines{conformal transformation} 
if $f$ is a conformal mapping that is onto $D'$. 
It follows that $f'(z) \neq 0$ for $z\in D$, and $f^{-1}:D' \to D$
is also a conformal transformation. 
We write $\univalent$ 
to denote the set of functions $f$ which are analytic and univalent 
in $\D$ satisfying the normalizing conditions $f(0)=0$ and $f'(0)=1$. 
 
If $D \subset \C$ with $0 \in D$,   we 
define the \defines{radius} (with respect to
 the origin) of $D$ to be $\rad(D)  := \sup\{|z| : z \in \bd D\}$,
 and the \defines{inradius} (with respect to the origin) of
 $D$ to be $\inrad(D) := \dist(0,\bd D) := \inf\{|z| : z \in \bd D\}$. 
The \defines{diameter} of $D$ is given by $\diam(D):=\sup\{|x-y| : x,y\in D\}$.
 If $D \subset \C $, then we say that a bounded $D$ is
a {\em Jordan domain}  if $\bd D$
 is a Jordan curve (i.e., homeomorphic to a circle).
  A Jordan domain is {\em nice}
if    the Jordan
 curve $\bd D$ can be expressed as a finite union of analytic curves. 
Note that Jordan domains  are simply connected.
For each $r>0$, let $\Dr{r}$ be the set of nice Jordan domains
 containing the origin of inradius $r$,
 and write $\scd := \bigcup_{r>0}\Dr{r}$. We also define 
$\scd^*$ to the be set of Jordan domains
containing the origin,   
and note that $\D \in \scd \subset \scd^*$. If $D$, $D' \in \scd^*$,
 let $\mathcal{T}(D,D')$ be the set of all $f : D \to D'$ that are  
conformal transformations of $D$ onto $D'$. The Riemann mapping theorem
 implies that $\mathcal{T}(D,D') \neq \emptyset$, and since $\bd D$, $\bd D'$
are Jordan, the Carath\'eodory extension theorem tells 
us that $f \in \mathcal{T}(D,D')$ can be extended to a 
homeomorphism of $\overline{D}$ onto $\overline{D'}$. 
We will use this fact repeatedly throughout,
 without explicit mention of it. For statements 
and details on these two theorems, consult~\cite[\S 1.5]{Duren}. 

A subset $A \subset \Z^2$  is \defines{connected} if every 
two points in $A$ can be connected by a nearest neighbour path
 staying in $A$.  We say that a finite subset $A$ is \defines{simply connected} 
if both $A$ and $\Z^2 \setminus A$ are connected.  There are three standard ways 
to define the ``boundary'' of a proper subset $A$ of $\Z^2$:
\begin{center}
\begin{itemize}
\item \textit{(outer) boundary:} $\;\bd A := \{y \in \Z^2 \setminus A: 
|y-x| = 1 \text{ for some } x \in A \}$;\\
\item \textit{inner boundary:} $\;\bd_i A :=\bd(\Z^2 \setminus A) =
 \{x\in A: |y-x|=1 \text{ for some } y \in \Z^2 \setminus A\}$;\\
\item \textit{edge boundary:} $\;\bd_e A := \{(x,y): x \in A,\, 
y \in \Z^2 \setminus A,\, |x-y| = 1\}$.
\end{itemize}
\end{center}

To each finite, connected
 $A \subset \Z^2$ we associate a domain  
$\tilde A \subset \C$ in the following way. 
For each edge $(x,y) \in \bd_e A$, considered as a 
line segment of length one, let $\ell_{x,y}$ be the perpendicular 
line segment of length one intersecting $(x,y)$ in the midpoint. 
 Let $\bd \tilde A$ denote the union of the line segments $\ell_{x,y}$,
 and let $\tilde A$ denote the domain with boundary $\bd \tilde A$ 
containing $A$. Observe that
\begin{equation}\label{squaredef}
\tilde A \cup \bd \tilde A = \bigcup_{x \in A}
 \Square_x \;\text{ where }\; \Square_x := x + 
\left(\,[-1/2,1/2] \times [-1/2,1/2]\,\right).
\end{equation}
That is, $\Square_x$ is the closed square of side 
length one centred at $x$ whose sides are parallel 
to the coordinate axes.  Also, note that $\tilde A$ is 
 simply connected if and only if $A$ is a simply
connected subset of $\Z^2$. We refer to $\tilde{A}$ as the
 \defines{``union of squares''} domain associated to $A$.

Let $\A$ denote the set of all finite simply connected 
subsets of $\Z^2$ containing the origin. If $A \in \A$, let 
$\inrad(A):=\min\{|z| : z \in \Z^2 \setminus A \}$ 
and $\rad(A):=\max\{|z| : z \in A\}$
denote the \defines{inradius} and \defines{radius} 
(with respect to the origin), respectively, of $A$, and 
define $\A^n$ to be the set of $A \in \A$ with $n \leq \inrad(A) \leq 2n$; 
thus $\A := \bigcup_{n\ge 0}\A^n$.  Note that if $A \in \A$ 
and  $0 \neq x \in \bd_i A$, then  the connected component of 
$A \setminus \{x\}$ containing the origin is simply connected. 
(This is not true if we do not assume $x \in \bd_i A$.)  
Similarly, by induction, if $A \in \A$, $0 \neq x_1 \in \bd_i A$,
 and $[x_1,x_2,\ldots,x_j]$ is a nearest neighbour path in $A \setminus \{0\}$, then
the connected component of  $A \setminus  \{x_1,\ldots,x_j\}$
 containing the origin is simply connected. 

Finally, if $A \in \A$ with associated domain $\tilde A \subset \C$, 
then we write
$f_A := f_{\tilde{A}}$ for the conformal transformation of $\tilde A$ onto 
the unit disk $\D$ with $f_A(0) =0$, $f_A'(0) > 0$.

\subsection{Green's functions on $\C$}\label{g-sect}

If $D$ is a domain whose boundary includes a curve, 
let $g_D(x,y)$ denote the Green's function for Brownian
motion. If $x \in D$, we can define $g_D(x, \cdot)$ 
as the unique harmonic function on $D\setminus \{x\}$, 
vanishing on $\bd D$ (in the sense that $g_D(x,y) \to 0$ as $y \to y_0$
for every regular $y_0 \in \bd D$), with 
\begin{equation}\label{greencharharm}
g_D(x,y) = -\log|x-y| + O(1) \;\text{ as }\; |x-y| \to 0. 
\end{equation}  In the case $D = \Disk$, we have
\begin{equation}\label{defngreen}
g_{\D}(x,y)=\log \left|\overline{y}x-1\right| - \log \left|y-x\right|.
\end{equation}
Note that $g_{\D}(0,x) = -\log |x|$, and $g_{\D}(x,y)=g_{\D}(y,x)$. 
The Green's function is a well-known example of a conformal invariant
(see, e.g.,~\cite[\S 1.8]{Duren}).

\begin{proposition}[Conformal Invariance of the Green's Function]\label{Green-conf}
Suppose that $f:D \rightarrow D'$ is a conformal transformation. If $x$, $y \in D$, then
$g_D(x,y) = g_{D'}(f(x),f(y))$.
\end{proposition}

 Hence, to determine $g_D$ for arbitrary $D \in \scd^*$, it is 
enough to find $f_D \in \mathcal{T}(D,\D)$. 
On the other hand, suppose $D \subset \C$ is a simply connected
 domain containing the origin
with Green's function $g_D(z,w)$.  Then we can write the unique
conformal transformation of $D$ onto $\Disk$ with $f_D(0) = 0,
f_D'(0) > 0$ as  
\begin{equation}\label{riemannmap}
f_D(x) = \exp\{-g_D(x) + i\theta_D(x)\},
\end{equation}
where $g_D(x) =g_D(0,x)$ and   $-g_D + i \theta_D$ is analytic in 
$D \setminus \{0\}$.

An equivalent formulation of the Green's function 
for $D \in \scd^*$ can be given in terms of Brownian motion. 
Using~(\ref{greencharharm})
 we have that 
$g_D(x,y) = \E^x[\log|B_{T_D}-y|] - \log |x-y|$
for distinct points $x$, $y \in D$ where $T_D := \inf \{t : B_t \not\in D \}$. 
In particular, if $0 \in D$, then 
\begin{equation}\label{1.3}  
g_D(x) = \E^x[\log|B_{T_D}|] - \log |x| \, \text{ for }\, x\in D.
\end{equation}
Additional details may be found in~\cite[Chapter 2]{LawlerSLE}\label{SLE1}.

If $A\in \A$, and if we 
let $g_A(x,y) := g_{\tilde A}(x,y)$ be the
 Green's function (for Brownian motion) in 
$\tilde A$, then by Proposition~\ref{Green-conf} 
and~(\ref{defngreen}) we have that  
\begin{equation}\label{sept9.1}
g_A(x,y)=g_{\D}(f_A(x), f_A(y))=\log
 \left|\frac{\overline{f_A(y)}f_A(x)-1}{f_A(y)-f_A(x)}\right|.
\end{equation}
If we write $\theta_A : =\theta_{\tilde{A}}$, then~(\ref{riemannmap}) implies that
\begin{equation}\label{oct5.1}  
f_A(x) = \exp\{-g_A(x) + i\theta_A(x)\}.
\end{equation}

\subsection{Green's functions on $\Z^2$}\label{greenZ}

Suppose that $S_n$ is a simple random walk on $\Z^2$ and
 $A$ is a proper subset
of  $\Z^2$. If $\tau_A := \min \{j \ge 0 : S_j \not\in A\}$,
 then  we let
\begin{equation}\label{SRWgreen}
G_A(x,y) := \E^x [\sum_{j=0}^{\tau_A -1} \mathbbm{1}_{\{S_j =y\}}] = 
\sum_{j=0}^{\infty} \P^x \{S_j=y, \tau_A > j \}
\end{equation}
denote the Green's function for random walk on $A$. Set $G_A(x) := G_A(x,0) = G_A(0,x)$. 
In analogy with the Brownian motion case (see~\cite[Proposition 1.6.3]{LawlerGreen}),
 we have
\begin{equation}\label{1.1}
G_A(x)=\E^x[a(S_{\tau_A})]-a(x) \, \text{ for }\, x\in A
\end{equation}
where $a$ is the potential kernel for simple random walk defined by
\[a(x)  :=  \sum_{j=0}^\infty \left[\P^0\{S_j=0\} - \P^x \{ S_j=0\}\right].\]
It is also known~\cite[Theorem 1.6.2]{LawlerGreen} that as $|x| \to \infty$,
\begin{equation}\label{1.2}
a(x) = \frac{2}{\pi}\log|x| + k_0 +o(|x|^{-3/2})
\end{equation}
where $k_0 := (2\varsigma + 3\ln 2)/\pi$ and $\varsigma$ is Euler's constant.  
The error above will suffice for our purposes, 
even though stronger results are known. 
The asymptotic expansion of $a(x)$ given in~\cite{FukU1} shows 
that the  error is $O(|x|^{-2})$.

\subsection{Consequences of the Koebe theorems}\label{koebesec}

We now recall some
standard
 results from the study of univalent functions. 
Proofs may be found in~\cite[Theorems 2.3, 2.4, 2.5, 2.6]{Duren}.

\begin{theorem}[Koebe One-Quarter Theorem]\label{koebethm}
If $f$ is a conformal mapping of the unit disk with $f(0) = 0$, 
then the image of $f$ contains the open disk of radius $|f'(0)|/4$ about the origin.
\end{theorem}

\begin{theorem}[Koebe Growth and Distortion Theorem]\label{growthdistortion}
If $f \in \univalent$ and $z \in \D$, then
$$\left|\frac{zf''(z)}{f'(z)} - \frac{2|z|^2}{1-|z|^2} \right| 
\le \frac{4|z|}{1-|z|^2}\, , \quad \frac{|z|}{(1+|z|)^2} \le |f(z)|
 \le \frac{|z|}{(1-|z|)^2} \, , $$
$$\frac{1-|z|}{(1+|z|)^3} \le |f'(z)| \le \frac{1+|z|}{(1-|z|)^3}.$$
\end{theorem}

A number of useful consequences may now be deduced.   

\begin{corollary}\label{koebecor1}
For each $0 < r< 1$, there is a constant $c_r$ such that if $f \in 
\univalent$ and $|z| \le r$, then $|f(z)-z| \le c_r |z|^2$.
\end{corollary}

\begin{proof}
If we combine the first estimate in Theorem~\ref{growthdistortion}
 with the estimate of $|f'(z)|$ in the third statement of that theorem, 
then we can obtain a uniform bound on $|f''(z)|$ over all $f \in \univalent$ 
and $|z| \le r$.
\end{proof}

Recall that $f_A :=f_{\tilde A} \in \mathcal{T}(\tilde{A},\D)$ is 
the unique conformal transformation of $\tilde A$ onto
$\Disk$ with $f_A(0)=0$, $f_A'(0)>0$.

\begin{corollary}\label{koebecor2}
If $A \in \A^n$, then $-\log f_A'(0) = \log n + O(1)$.
\end{corollary}

\begin{proof}
Using the Koebe one-quarter theorem and the Schwarz lemma, we see that 
if $F: \D \to D$ is a conformal transformation with $F(0) = 0$, then
\[     |F'(0)|/4 \leq \inrad(D) \leq |F'(0)|. \]
By definition, since $A \in \A^n$, 
we have that $n-1 \le \inrad(\tilde{A}) \le 2n+1$. 
Therefore, if $F_{A} :=f^{-1}_{A}$, then
\begin{equation*}
          n-1 \leq F_A'(0)  \leq 4(n+1). \qedhere
\end{equation*}
\end{proof}

Along with Corollary~\ref{koebecor1}, the growth and distortion theorem yields the following.

\begin{corollary}\label{koebecor3}
If $A \in \A^n$ and $|x| \leq n/16$, then $f_A(x) = x f_A'(0) + |x|^2 O(n^{-2})$,
and
\begin{equation}
g_A(x) + \log|x| = - \log f_A'(0) + |x|O(n^{-1}).
\end{equation}
\end{corollary}

\begin{proof}
For $z \in \D$, let $F_A(z) := f_A(nz)/(nf'_A(0))$. 
 Then $F_A \in \univalent$, so Corollary~\ref{koebecor1}
 with $r=1/16$ gives $|F_A(z) - z | \le C |z|^2$. Thus, 
if $z = x/n$, $|f_A(x) - x f_A'(0)| \le C f'_A(0)|x|^2 n^{-1}$.
 By the previous corollary, $f'_A(0) = O(n^{-1})$, so the first
assertion follows.
The second result follows
from $|f_A(x)| = \exp\{-g_A(x)\}$.
\end{proof}

We remark that this corollary implies  
$\lim_{|x|\to 0} (g_A(x) + \log|x| ) = - \log f_A'(0)$
which shows the size of
 the error term in~(\ref{greencharharm}). 

\subsection{Beurling estimates and related results}\label{beurlingsec}

Throughout this subsection, suppose that $A \in \A^n$ with
 associated ``union of squares'' domain $\tilde{A} \subset \C$, 
and write $T_A := T_{\tilde{A}} := \inf \{t : B_t \not\in \tilde{A}\}$. 
From the Beurling projection theorem~\cite[Theorem~(V.4.1)]{Bass} the 
following may be derived.

\begin{theorem}\label{theorembeurlingestimate}\label{corollarybeurlingestimate1}
There is a constant $c < \infty$ such that if 
$\gamma:[a,b] \to \C$ is a curve with $|\gamma(a)| =r$, 
$|\gamma(b)|=R$, $0 < r < R <\infty$, $\gamma(a,b) \subset
 \D_R:=\{z \in \C : |z| < R\}$, and $|x|\le r$, then
\begin{equation}\label{beurcor1}
\P^{x}\{B[0,T_{\D_R}] \cap \gamma[a,b] = \emptyset \} \le c \: (r/R)^{1/2}.
\end{equation}
\end{theorem}

\begin{corollary}[Beurling Estimate]\label{corollarybeurlingestimate2}
There is a constant $c<\infty$ such that if $x \in \tilde{A}$, then for all $r>0$,
\begin{equation}\label{beurcor2}
\P^{x}\{|B_{T_A} - x | > r \: \dist(x,\bd \tilde A)\} \le c \:r^{-1/2}.
\end{equation}
\end{corollary}

\begin{proof}
Without loss of generality, we may assume 
by Brownian scaling that $x=0$ and $\inrad(\tilde{A}) =:
 d \in [1/2, 1]$.  If $\rad(\tilde{A}) \le r$, then this 
estimate is trivial.  If not, then there is a curve in $\bd \tilde{A}$ 
from the circle of radius $d$ to the circle of radius $r$, and the Beurling estimate~(\ref{beurcor2}) 
follows from~(\ref{beurcor1}).
\end{proof}

In particular,  if $|x| > n/2$, the probability
 starting at $x$ of reaching $\D_{n/2}$ before leaving $\tilde A$ i
s bounded above by $cn^{-1/2} \dist(x,\bd\tilde A)^{1/2}$. From the
 Koebe one-quarter theorem it easily follows that $g_A(x) \leq c$ for 
$ |x| \geq n/4$; hence we get
\begin{equation}\label{gbound}
g_A(x) \leq c \: n^{-1/2}\: \dist(x,\bd\tilde A)^{1/2}, \;\;\;\; A \in 
\A^n, \;\;\; |x| \geq n/4.
\end{equation}

Recall from~(\ref{oct5.1}) that $f_A(x) = \exp\{-g_A(x) + 
i\theta_A(x)\}$ for $x \in \tilde{A}$. Hence, if $x \in \bd_i A$,  
then $g_A(x) \leq c n^{-1/2}$, so that $f_A(x) = \exp\{i\theta_A(x)\} +
 O(n^{-1/2})$. If $z \in \bd A$, then since $f_A(z)$ is not defined,
 we let $\theta_A(z)$ be the average of $\theta_A(x)$ over all $x \in A$ 
(for which $f_A(x)$ \emph{is} defined) with $|x-z| = 1$.  The Beurling estimate 
and a simple Harnack principle show that
\begin{equation}\label{theta_estimate}
\theta_A(z) = \theta_A(x) + O(n^{-1/2}),  \;\;\;\;  (x,z) \in \bd_eA.
\end{equation}

There are analogous Beurling-type results in the discrete case; the 
following is a corollary of~\cite[Theorem~2.5.2]{LawlerGreen}. 
Let $\tau_A := \min \{j \ge 0 : S_j \not\in A\}$. 

\begin{corollary}[Discrete Beurling Estimate]\label{discretebeurlingest}
There is a constant $c < \infty$ such that if $r > 0$, then 
$\P^{x}\{|S_{\tau_A} - x | > r \: \dist(x,\bd A)\} \le c \:r^{-1/2}$.
\end{corollary}

In particular,  if $|x| > n/2$, the probability 
starting at $x$ of reaching $\D_{n/2}$ before leaving $A$ is 
bounded above by $cn^{-1/2} \dist(x,\bd A)^{1/2}$. It is easy to show 
that $G_A(x) \leq c$ for $ |x| \geq n/4$; hence in this case we get
\begin{equation}\label{Gbound}
G_A(x) \leq c \: n^{-1/2}\: \dist(x,\bd A)^{1/2}, \;\;\;\; A \in \A^n, \;\;\; |x| \geq n/4.
\end{equation}
Specifically, if $x \in \bd_i A$,  then $G_A(x) \leq c n^{-1/2}$.

If $A \in \A$ and $0 \neq x \in \bd_i A$, then since 
$G_A(0) = G_{A \setminus \{x\}}(0) + \P\{\tau_A > \tau_{A\setminus
 \{x\}}\}G_A(x)$ it follows that
$$G_A(0) = G_{A \setminus \{x\}}(0) + \frac{G_A(x)^2}{G_A(x,x)}.$$ 
We can replace $A \setminus \{x\}$ in the above formula with the 
connected component of $A \setminus \{x\}$ containing the origin. 
In particular, if $A \in \A^n$ and $x \in \bd_i A$, then we 
conclude that $G_A(0) - G_{A \setminus \{x\}}(0) \leq G_A(x)^2 \leq c\;n^{-1}$.

\subsection{Excursion Poisson kernel}

Let $D$ be a domain in $\C$.  We
say that a connected
$\Gamma$ is an {\em (open) analytic arc} of $\p D$ if there
is a domain $D' \subset \C$ that is symmetric
about the real axis  and a
conformal map $f: D' \rightarrow f(D')$ such that
$f(D' \cap \R) = \Gamma$ and $f(D' \cap \Half) = f(D') \cap D$, where
$\Half$ denotes the upper half plane.  We
say that $\p D$ is {\em locally analytic} at $x \in \p D$, if there
exists an analytic arc of $\p D$ containing $x$. 
 If $\Gamma$ is an analytic 
arc of $\p D$ and $x,y \in \Gamma$, we write $x < y$ if $f^{-1}(x)
< f^{-1} (y)$.  We let $\mathbf{n}_x := \mathbf{n}_{x,D}$ be the unit normal
at $x$ pointing into $D$.

If $x \in D$ and
$\p D$ is locally analytic at  $y \in \bd D$, 
then both \defines{harmonic measure} $\P^{x} \{B_{T_D} \in \d y \}$, and its density with 
respect to arc length, the \defines{Poisson kernel} $H_D(x,y)$, are well-defined. Also,
 recall that for fixed $z \in D$, the function $y \mapsto H_D(z,y)$ is continuous 
in $y$, and that for fixed $y\in \bd D$, the function $z \mapsto H_D(z,y)$ is 
harmonic in $z$. If $x \in D$ and $\Gamma \subset \bd D$ is an analytic  
  arc, then write
\begin{equation}\label{PK-int-new}
H_{D}(x,\Gamma) := \int_{\Gamma} H_{D}(x,y) \, |\d y|.
\end{equation}
The Riemann mapping theorem and L\'evy's theorem on the conformal
 invariance of Brownian motion~\cite{Bass} allow us to describe the 
behaviour of the Poisson kernel under a conformal transformation. 

\begin{proposition}\label{PKconf}
If $f: D \to D'$ is a conformal transformation,
$x \in D$,  $\p D$ is locally analytic at $y \in \bd D$,
and $\p D'$ is locally analytic at $f(y)$,
then  
$\P^{x} \{B_{T_D} \in \d y \} = \P^{f(x)} \{{B'}_{T_{D'}} \in f(\d y) \}$
where $B'$ is a (time-change of) Brownian motion. Equivalently,
\begin{equation}\label{conf-inv-PK}
H_D(x,y) = |f'(y)| \, H_{D'}(f(x), f(y)).
\end{equation}
\end{proposition}
For each $\epsilon > 0$, let
$\mu_{x,\epsilon,D}$ denote the probability measure on paths obtained
by starting a Brownian motion at $x + \epsilon \mathbf{n}_x$ and stopping
the path when it reaches $\p D$.  The excursion measure of $D$ at
$x$ is defined by
\[    \excur_D(x) := \lim_{\epsilon \rightarrow 0+} \frac{1}{\epsilon} \,
   \mu_{x,\epsilon,D}. \]
Excursion measure on $\Gamma$ is defined by
\[    \excur_D(\Gamma) := \int_\Gamma  \excur_D(x) \, |dx| . \]
If $\Upsilon$ is another analytic arc on $\p D$, we define
$\excur_D(\Gamma,\Upsilon)$ to be the excursion measure $\excur_D
(\Gamma)$ restricted to curves that end at $\Upsilon$ and whose
endpoints are different.   The {\em excursion boundary measure} is defined by
\[    H_{\bd D}(\Gamma,\Upsilon) := |\excur_D(\Gamma,\Upsilon)|, \]
where $|\cdot|$ denotes total mass.  We can write
\[    H_{\bd D}(\Gamma,\Upsilon) = \int_\Upsilon \int_\Gamma 
    H_{\bd D}(x,y) \, |dx| \, |dy| , \]
where $H_{\bd D}(x,y)$, $x \neq y$, denotes the {\em excursion Poisson kernel} given by
\[    H_{\bd D}
(x,y) := \lim_{\epsilon \rightarrow 0+}
                  \frac{1}{\epsilon} \, H_{  D}(x + \mathbf{n}_x,y)
                   = \lim_{\epsilon \rightarrow 0+}
                  \frac{1}{\epsilon} \, H_{  D}( y + \mathbf{n}_y,x). \]
   We can also write
\[    \excur_D(\Gamma,\Upsilon) =  \int_\Upsilon \int_\Gamma
    H_{\bd D}(x,y) \, \excur_D^\#(z,w) \, |dx| \, |dy| , \]
where $\excur_D^\#(z,w) := \excur_D(z,w)/H_{\bd D}(x,y)$ is the
excursion measure between $x$ and $y$ normalized to be a probability
measure.  We will consider $\excur_D$ and $\excur_D^\#$ as measures on
curves modulo reparametrization.  Conformal invariance of
complex Brownian motion implies that $\excur_D^\#$ is
conformally invariant, i.e., if $f: D \rightarrow D'$ is a conformal
transformation, then
\[   f \circ \excur_D^\#(x,y) =  \excur_{D'}^\#(f(x),f(y)). \]

\begin{proposition}\label{EPKconf}
Suppose   $f: D \to D'$ is a conformal 
transformation and $x,y$ are distinct points
on $\p D$.  Suppose that $\bd D$ is locally analytic
at $x,y$ and $\bd D'$ is locally analytic
at $f(x), f(y)$.
Then
$H_{\bd D}(x,y) = |f'(x)| \,|f'(y)| \, H_{\bd D'}(f(x), f(y))$.
\end{proposition}

\begin{proof}
By definition, $H_{\bd D}(x,y) := \lim_{\e \to 0+} \e^{-1} \, H_{D}(x+\e\n_x,y)$. Therefore,
\begin{align*}
H_{\bd D}(x,y)&=\lim_{\e \to 0+} \e^{-1} \, |f'(y)| \, H_{D'}(f(x+\e\n_x),f(y)) 
\;\text{ (Proposition~\ref{PKconf})}\\
&=|f'(y)| \, \lim_{\e \to 0+} \e^{-1} \, H_{D'}(f(x)+\e f'(x)\n_x +o(\e), f(y)) 
 \\ 
&=|f'(x)| \, |f'(y)| \,\lim_{\e_1 \to 0+} \e_1^{-1} \, H_{D'}(f(x)+ \e_1\n_{f(x)}+o(\e_1),f(y))\\
&=|f'(x)| \, |f'(y)| \,H_{\bd D'}(f(x), f(y)) 
\end{align*}
where we have written $\e_1:=\e|f'(x)|$, and have noted that $\n_{f(x)} = f'(x)\n_x/|f'(x)|$ and $|\n_{f(x)}|=|\n_x|=1$.
\end{proof}

\begin{example}\label{ex1}
If $x = e^{i \theta}, y = e^{i \theta'}$, $\theta \neq \theta'$, then 
\[ H_{\bd \Disk}(x,y) = \frac 1 \pi \, \frac 1 {|y-x|^2}
  = \frac 1 {2\pi} \, \frac{1}{1 - \cos(\theta' - \theta)}. \]
\end{example}


\begin{example}\label{ex2}
If $\rect_L = \{z: 0 < \Re(z) < L, 0 < \Im(z) < \pi\}$, then
separation of variables can be used to show that for $0<r<L$, $0<q,q'<\pi$,
\[   H_{ \rect_L}(r + iq,L + iq')
              = \frac 2 \pi \sum_{n=1}^\infty 
     \frac{\sinh(nr) \, \sin(nq) \, \sin(nq')}{\sinh (nL)}, \]
\begin{equation}  \label{dec8.1}
  H_{\bd \rect_L}(iq, L + i q') = \frac 2 \pi \sum_{n=1}^\infty
    \frac{n \, \sin(nq) \, \sin(nq')}{\sinh (nL)} . 
\end{equation}
\end{example}

\begin{example}\label{ex3}
In the Brownian excursion case, there is an exact form of Theorem~\ref{rwestimate}.
Suppose that $D \in \scd$, and $x$, $y \in \bd D$ with $\bd D$ locally analytic at $x$ and $y$. 
Proposition~\ref{EPKconf} and Example~\ref{ex1} imply that 
$2\pi H_{\bd D}(x,y) = |f'(x)| \,|f'(y)| \, (1 - \cos(\theta_D(x) - \theta_D(y)))^{-1}$ where $f \in \mathcal{T}(D, \D)$ with $f(0)=0$.
However, Proposition~\ref{PKconf} combined with the fact that harmonic measure from 0 in $\D$ is uniform on $\bd \D$, 
implies $|f'(x)|^{-1} \, H_{D}(0,x) = |f'(y)|^{-1} \,H_{D}(0,y)=(2\pi)^{-1}$. Hence,   
we conclude
\begin{equation*}
H_{\bd D}(x,y) =\frac{2\pi \, H_{D}(0,x) \, H_{D}(0,y)}{1- \cos(\theta_D(x)-\theta_D(y))}.
\end{equation*}
\end{example}
 
We finish this subsection by
stating a formula relating the excursion 
Poisson kernel and the Green's function for Brownian 
motion:
\begin{equation}\label{H-g-analogue}
H_{\bd D}(x,y) =\lim_{\e \to 0+} \frac{g_{D}(x+\e\n_x,y+\e\n_y)}
{2\pi\,\e^2}= \frac{\bd}{\bd \n_x} \frac{\bd}{\bd \n_y} g_D(x,y).
\end{equation}
For the simply connected case, which is all that we will
use, this follows from a straightforward computation.  We omit the
proof in the general case. 

\subsection{Excursion Poisson kernel determinant}

\begin{definition}\label{hitmatrixD-b}
If $D$ is a domain, and  $x^1,\ldots,x^k,y^1,\ldots,y^k$
are distinct boundary points at which $\bd D$ is locally
analytic,
let $\mathbf{H}_{\bd D}
(\mathbf{x},\mathbf{y}) = 
[H_{\bd D}(x^i, y^j)]_{1\le i,j \le k}$ denote the 
\defines{$k \times k$ hitting matrix}
\begin{equation*}
\mathbf{H}_{\bd D}(\mathbf{x},\mathbf{y}) := 
\begin{bmatrix}
H_{\bd D}(x^1, y^1) &\cdots & H_{\bd D}(x^1, y^k) \\
\vdots      &\ddots &\vdots       \\
H_{\bd D}(x^k, y^1) &\cdots & H_{\bd D}(x^k, y^k)
\end{bmatrix}
\end{equation*}
where $H_{\bd D}(x^i,y^j)$ is the excursion Poisson kernel, and $\mathbf{x} =
(x^1,\ldots,x^k)$, $\mathbf{y} = (y^1,\ldots,y^k)$.
\end{definition}

A straightforward extension of Proposition~\ref{EPKconf} is that the
 determinant of the hitting matrix of excursion Poisson kernels is conformally covariant.

\begin{proposition}\label{poisson_det_lemma-b}
If $f: D \rightarrow D'$ is a conformal
transformation, $x^1,\ldots,x^k$,$y^1,\ldots,y^k$ are distinct points at which $\bd D$ is locally analytic, and
$\bd D'$ is locally analytic at $f(x^1),\ldots,f(x^k)$, $f(y^1),\ldots,f(y^k)$, then
$$\det\,\mathbf{H}_{\bd D}(\mathbf{x},\mathbf{y}) = \left(\prod_{j=1}^{k}|f'(x^j)| \, 
|f'(y^j)| \right) \det\, [H_{\bd D'}(f(x^i), f(y^{j}))]_{1 \le i,j \le k}.$$ 
\end{proposition}

\begin{proof}
By the definition of determinant,
\begin{align*}
&\det\, \mathbf{H}_{\bd D}(\mathbf{x},\mathbf{y})
  =   \sum_{\sigma} (\sgn \, \sigma)  \: 
H_{\bd D}(x^1, y^{\sigma(1)}) \cdots H_{\bd D}(x^k, y^{\sigma(k)})\\
&= \left(\prod_{j=1}^{k}|f'(x^j)| \, |f'(y^j)| \right) \sum_{\sigma} 
(\sgn \, \sigma) \: H_{\bd D'}(f(x^1), f(y^{\sigma(1)})) \cdots H_{\bd D'}(f(x^k),
 f(y^{\sigma(k)}))\\
&= \left(\prod_{j=1}^{k}|f'(x^j)| \, |f'(y^j)| \right) \det\,
 [H_{\bd D'}(f(x^i), f(y^{j}))]_{1 \le i,j \le k}
\end{align*}
where the sum is over all permutations $\sigma$ of $\{1,\ldots,k\}$,
and $\sgn \, \sigma$ denotes the sign of the permutation.
\end{proof} 

It follows from this proposition that $\excdetfrac_D(x^1,\ldots,x^k,y^k,
\ldots,y^1)$ given by~(\ref{conditionalFominBM}) is a conformal invariant, and hence it can be defined
for any Jordan domain by conformal invariance. 

\subsection{Proof of Proposition~\ref{propdec13.1}}

By applying  Proposition~\ref{poisson_det_lemma-b} and judiciously choosing $D=\rect_L$, we now give
the asymptotics as points gets close, and complete the proof of  Proposition~\ref{propdec13.1}.

\begin{proof}[Proof of Proposition~\ref{propdec13.1}]  
Let $\mathbf{q}=(q_1,\ldots,q_k)$, $\mathbf{q}'=(q_1',\ldots,q_k')$,
\[  \vec u_j = \left[\begin{array}{c} \sin(jq_1) \\
  \sin (j q_2) \\ \vdots \\ \sin(j q_k) \end{array} \right],
\;\;\;\;\;\;  \vec v_j = \left[\begin{array}{c} \sin(jq_1') \\
  \sin (j q_2') \\ \vdots \\ \sin(j q_k') \end{array} \right].\] 
Using~(\ref{dec8.1}), we see that 
$(\pi/2)^k\,\det\,\mathbf{H}_{\bd \rect_L}(i\mathbf{q},L + i\mathbf{q}')$ can be written
as
\[  
\det\left[\sum_{j=1}^\infty \frac{ j\,\sin(j q_1)}{\sinh(jL)} \, \vec v_j,
      \sum_{j=1}^\infty \frac{j\,\sin(j q_2)}{\sinh(jL)} \, \vec v_j,
\ldots, \sum_{j=1}^\infty \frac{j\, \sin(j q_k)}{\sinh(jL)}
 \, \vec v_j\right]_{1\le j \le k}. \]
By multilinearity of the determinant, we can write the determinant
above as
\[    \sum_{j_1, \ldots,j_k} 
          \frac{  (j_1 \cdots j_k)  \; \sin(j_1 q_1) \cdots 
   \, \sin(j_k q_k)}{\sinh(j_1L) \cdots
  \sinh(j_kL) } \, \det[\vec v_{j_1},\ldots, \vec v_{j_k}] . \]
The determinants  in the last sum
equal zero   if the  indices are not distinct. Also
it is not difficult to show that
\[  \sum_{j_1+ \cdots + j_k \geq R}
                \frac{   j_1 \cdots j_k     } {\sinh(j_1L) \cdots
  \sinh(j_kL) } \leq  C(k,R) \, e^{-RL}. \]
Hence, except for   an error of $O_k(e^{-(k^2 + k + 2)L/2})$, we see that
 $(\pi/2)^k \,
\det\,\mathbf{H}_{\bd \rect_L}(i\mathbf{q},L + i\mathbf{q}')$ equals
\begin{equation}  \label{dec2.1}
  k! \,  \sum_{\sigma} \frac{
 \sin ( {\sigma(1)}q_1) \cdots \sin( {\sigma(k)}q_k)  }
  {\sinh(L) \, \sinh(2L) \cdots \sinh(kL)}   
   \,   \det[\vec v_{\sigma(1)},\ldots, \vec v_{ {\sigma(k)}}], 
\end{equation}
where the sum is over all permutations $\sigma$ of $\{1,\ldots,k\}$.
But \[\det[\vec v_{\sigma(1)},\ldots, \vec v_{ {\sigma(k)}}] =
(\sgn \, \sigma) \, \det[\vec v_1,\ldots,\vec v_k] .\]
Hence~(\ref{dec2.1}) equals
\[    \frac{  k!  
\;  \det[\vec u_1,\ldots,\vec u_k] \; \det[\vec v_1,\ldots,\vec v_k]}
   {\sinh(L) \, \sinh(2L) \cdots \sinh(kL)} ,  \]
which up to 
an error of $O_k(e^{-(k^2 + k + 2)L/2})$ equals
\[    2^k \, k! \; e^{-k(k+1)L/2}
\; \det[\vec u_1,\ldots,\vec u_k] \; \det[\vec v_1,\ldots,\vec v_k] . \]
To finish the proof, note that from~(\ref{dec8.1}), we can also write
\[   H_{\bd \rect_L}(iq,L + iq') = 
  \frac 4 \pi \; e^{-L} \;  \sin q \; \sin q'   \;
   [1 + O(e^{-L})], \]
so that 
\begin{equation*}
  (\pi/2)^k \, \prod_{j=1}^k H_{\bd \rect_L}(iq_j, L + i q_j')
            = 2^k \, e^{-kL} \,
\left(\prod_{j=1}^k \sin q_j \, \sin q_j'\right) \;[1 + O_k(e^{-L})]. \qedhere
\end{equation*}
\end{proof}

\subsection{The discrete excursion Poisson kernel}\label{discrete_epk_sect}

\begin{definition}\label{defn-discretePK}
Suppose that $A$ is a proper subset of $\Z^2$,
 and let $\tau_A := \min\{j >0 : S_j \in \bd A \}$.  
For $x \in A$ and $y \in \bd A$, define the \defines{discrete Poisson kernel} 
(or \defines{hitting probability}) to be 
$h_A(x,y) := \P^{x}\{S_{\tau_A} = y\}$.
\end{definition}

We define the discrete analogue of the excursion Poisson kernel 
to be the probability that a random walk starting at $x \in \bd A$ takes 
its first step into $A$, and then exits $A$ at $y$.

\begin{definition}\label{discreteEPKdefn}
Suppose that $A$ is a proper subset of $\Z^2$,
 and let $\tau_A := \min\{j >0 : S_j \in \bd A \}$.  
For $x$, $y \in \bd A$, define the \defines{discrete excursion Poisson kernel}
 to be $h_{\bd A}(x,y) := \P^{x}\{S_{\tau_A} = y , S_1 \in A \}$.
\end{definition}

Note that
\[
h_{\bd A}(x,y)  = \frac{1}{4} \sum_{(z,x) \in \bd_e A} h_A(z,y).
\]
Also, a last-exit decomposition gives
\[
h_A(x,y) = \frac{1}{4}\sum_{(z,y) \in \bd_e A} G_A(x,z)
\]
where $G_A$ is the Green's function for simple random walk on $A$ as 
in~(\ref{SRWgreen}). 
Hence,
\begin{equation}  \label{dec15.1}
h_{\p A}(x,y) = \frac 1 {16} \sum_{ \; (z,y) \in \bd_e A\; }
          \sum_{\; (w,x) \in \bd_e A\; } G_A(w,z), 
\end{equation}
which is a discrete analogue of~(\ref{H-g-analogue}).


\section{Proofs of the Main Results}\label{Sect3}

\noindent After briefly reviewing strong approximation and establishing some ancillary results in the first subsection, we 
devote the second subsection to the proof of Theorem~\ref{greentheoremB} which relates
the Green's function for Brownian motion to the Green's function for simple 
random walk in certain domains.
 While this result  holds in general, it is  most useful 
for points away from the boundary.  In the final two subsections, we prove Theorem~\ref{rwestimate}
by obtaining  better estimates for the case 
of points near the boundary, provided they are not too close to each other. 
  Throughout this section, suppose that $A\in \A^n$ with 
associated ``union of squares'' domain $\tilde{A} \in \scd$. 
As in Section~\ref{scs}, let $f_A \in \mathcal{T}(\tilde{A},\D)$ with $f_A(0)=0$, 
$f'_A(0)>0$, and recall from~(\ref{oct5.1}) 
that $f_A(x) = \exp\{-g_A(x) + i\theta_A(x)\}$, where $g_A$ is 
the Green's function for Brownian motion in $\tilde{A}$. 

\subsection{Strong approximation}

In order to establish our Green's function estimates, we will need to establish a
 strong approximation result making use of the  theorem of Koml\'os, Major, 
and Tusn\'ady~\cite{KomMT1, KomMT2}.  
 Since we are concerned exclusively with complex Brownian 
motion and simple random walk
 on $\Z^2$, the results noted in~\cite{Aue1} suffice. 

\begin{theorem}[Koml\'os-Major-Tusn\'ady]\label{KMTthm}
There exists     $c <\infty$  and   a probability
 space $(\Omega, \F, \P)$ on which are    defined a two-dimensional Brownian 
motion $B$ and a two-dimensional simple random walk $S$ with $B_0=S_0$, such
 that for all $\lambda >0$ and each $n \in \N$,
\[\P \left\{\max_{10\leq t   \le n} \left|\frac{1}{\sqrt{2}}\,B_t - S_t\right| >
 c \,(\lambda + 1) \, \log n  \right\} < c \, n^{-\lambda}.\]
\end{theorem}

Here $S_t$ is defined for noninteger $t$ by linear interpolation.
The one-dimensional proof may 
be found in~\cite{KomMT2}, and the immediate extension 
to two dimensions is written down in~\cite[Lemma~3]{Aue1}.  
For our purposes, we will need to consider the maximum up to a \emph{random} 
time, not just a fixed time.  The following  strong approximation will suffice.
 Our choice of $n^8$ is arbitrary  and will turn out to be good enough.  

\begin{corollary}[Strong Approximation]\label{strapprox}
There  exist $C< \infty$
and  a probability space $(\Omega, \F, \P)$ on which are  
defined a two-dimensional Brownian motion $B$ and a
 two-dimensional simple random walk $S$ with $B_0=S_0$  such that
$$\P \{\max_{0 \le t \le \sigma_n} 
  |\frac 1 {\sqrt 2}B_t - S_t| > C \log n \} = O(n^{-10}),$$
where $\sigma^{1}_n := \inf\{t : |S_t-S_0|
 \ge n^8\}$, $\sigma^{2}_n := \inf\{t : |B_t-B_0| \ge n^8\}$,
 and $\sigma_n := \sigma^{1}_n \vee \sigma^{2}_n$.
\end{corollary}

\begin{proof}
If we choose $\rho > 0$ such that $\P\{|B_1| \geq 2 \}
\geq \rho$ nd $\P\{|S_n| \geq 2 \sqrt n \}\geq 
\rho$ for all $n$ sufficiently large, then iteration shows
that $\P\{\sigma_n > n^{36}\} \leq (1-\rho)^{n^{20}}
=o(n^{-10}).$   
Suppose that $\lambda = 5/18$, and let $C = 23c/18$ where
$c$ is 
as in Theorem~\ref{KMTthm}.  Then, by using that theorem,
\[
 \P \{\max_{0 \le t \le \sigma_n} 
  |\frac 1 {\sqrt 2}B_t - S_t| > C \log n \}
  \leq  
  \P \{\max_{0 \le t \le n^{36}} 
  |\frac 1 {\sqrt 2}B_t - S_t| > C \log n \}   
  + o(n^{-10})  
 =  O(n^{-10}),
\]
and the proof is complete.           
\end{proof}

\begin{proposition}\label{sept23.thm1.nowprop}
There exists a constant $c$ such that for every $n$, 
a Brownian motion $B$ and a simple random walk $S$ can be defined on
the same probability space so that if $A \in \A^n$, $1 < r \leq n^{20}$, 
and $x \in A$ with $|x| \leq n^3$, then
$$\P^x\{|B_{T_A} - S_{\tau_A}| \geq c r \log n \} \leq  c  r^{-1/2}.$$
\end{proposition}

\begin{proof}  
For any given $n$, let $ \hat
B$ and $S$ be defined as in Corollary~\ref{strapprox} 
above, and let $C$ be the constant in that corollary. 
Define $T'_A := \inf\{t \geq 0 : \dist(\hat B_t/\sqrt 2,\bd \tilde A)
 \leq 2 C \log n \}$, $\tau'_A := \inf\{t \geq 0 : \dist(S_t,\bd A)
 \leq 2 C \log n \}$, $\hat T_A = \inf\{t \geq 0:
\hat B_t / \sqrt 2 \in \bd A\} .$ and consider the events
\[  V_1 := \{\sup_{0 \leq t \leq  \sigma_n}:
|\frac 1{\sqrt 2} {\hat B_t}  - S_t| > C  \log n \}, \]\[
V_2 := \{\sup_{T'_A \leq t \leq \hat T_A} |
  \hat B_t -  \hat B_{T'_A} | \geq r  \log n \},\]
$$\text{and}\; V_3 := \{\sup_{\tau'_A \leq t \leq \tau_A} 
|S_t - S_{\tau'_A} | \geq r \log n\}.$$ 
By the Beurling estimates (Corollaries~\ref{corollarybeurlingestimate2}
 and~\ref{discretebeurlingest}), and the strong Markov property
(applied separately to the random walk and the Brownian motion), it follows 
that $\P(V_2 \cup V_3) = O(r^{-1/2})$. From Corollary~\ref{strapprox},
 $\P(V_1) =O(n^{-10}) = O(r^{-1/2})$. Therefore, $\P(V_1 \cup V_2 \cup V_3) =
 O(r^{-1/2})$.
 Note that $|(\hat B_{\hat T_A}/\sqrt 2) - S_{\tau_A}| 
\leq (r + 2C) \log n$ on the complement of $V_1 \cup V_2 \cup V_3$.
Also, if $B_t = \hat B_{2t}/\sqrt 2$, then $B_t$ is 
a standard Brownian motion and  $B_{T_A} = \hat B_{\hat T_A}
/ \sqrt 2$.
\end{proof}

\subsection{Proof of Theorem~\ref{greentheoremB}}\label{G-and-g-sec}

\begin{lemma}\label{sept9.lemma1}
There exists a $c$ such that if $A \in \A^n$ and $|x| \leq n^{2}$,  then
$$| \; \E^
x[\log |B_{T_A}|] - \E^x[\log |S_{\tau_A}|] \; |\leq c \; n^{-1/3}  \log n.$$
\end{lemma}
 
\begin{proof}
For any $n$, let $B$ and $S$ be as in Proposition~\ref{sept23.thm1.nowprop}, 
and let $V = V(n)$ be the event that $|B_{T_A} - S_{\tau_A}| \leq n^{2/3}
 \log n$.  By that proposition, $\P(V) \geq 1 - c n^{-1/3}$. Since $\inrad(A) \geq n$, 
we know that on the event $V$,
$$| \; \log| B_{T_A}| - \log |S_{\tau_A}|\; |\leq c' \;  n^{-1/3} \; \log n .$$
Note that 
$\E^x[\; \log| B_{T_A}|\; \mathbbm{1}_{V^c} \;\mathbbm{1}_{\{|B_{T_A}|
 \leq n^5\} } \; ] \leq c \; \log n \; \P(V^c) \leq c \; n^{-1/3} \; \log n$.
 Using the Beurling estimates it is easy to see that
$$\E^x[ \; \log| B_{T_A}| \; \mathbbm{1}_{ \{|B_{T_A}| \geq n^5\}} \; ] =
 O (n^{-1/3} \log n),$$
and similarly for $\log |S_{\tau_A}|$ in the last two estimates.  Hence,
\begin{equation*}
\E^x[\; (\; \log| B_{T_A}| + \log |S_{\tau_A}|\; ) \; \mathbbm{1}_{V^c} \; ]
 \leq c \; n^{-1/3}\;  \log n. \qedhere
\end{equation*}
\end{proof}

\begin{proof}[Proof of Theorem~\ref{greentheoremB}]
First suppose $x \neq 0$.
Recall from~(\ref{1.3}) that
 $g_A(x) = \E^x[\log|B_{T_A}|] - \log |x|$ and from~(\ref{1.1}) 
that $G_A(x) = \E^x[a(S_{\tau_A})] - a(x)$ 
with $a(x)$ as in~(\ref{1.2}). If $|x| \le n^2$, then~(\ref{mar29.eq1}) 
follows from Lemma~\ref{sept9.lemma1}, and if $|x| \ge n^2$, 
then~(\ref{mar29.eq1}) follows directly from the bounds on $g_A(x)$ 
and $G_A(x)$ in~(\ref{gbound}) and~(\ref{Gbound}), respectively.

If $x = 0$, we can use the relation
\[  G_A(0) = 1 + G_A(e) = a(1) + G_A(e) , \;\;\;\; |e| = 1, \]
and $|f_A(e)| =  |f_A'(0)| + O(n^{-2})$.
\end{proof}
 
For any $A \in \A^n$, let $A^{*,n}$ be the set
\begin{equation}\label{stardomain}
A^{*,n} := \{ x \in A: g_A(x) \geq n^{-1/16} \}.
\end{equation}
The choice of $1/16$ for the exponent is somewhat arbitrary, and slightly
 better estimates might be obtained by choosing a different exponent. However,
since we do not expect the error estimate derived here to be optimal, we will
 just make this definition.
 
\begin{corollary}\label{sept9.prop2.nowcor}
If $A \in \A^n$, and  $x\in A^{*,n}$, $y \in A$, then
\begin{equation}\label{mar29.eq2}
G_A(x,y) = \frac{2}{\pi} \, g_A(x,y) + k_{y-x} + O(n^{-7/24} \log n).
\end{equation}
\end{corollary}

\begin{proof} 
From the Beurling estimates,~(\ref{gbound}), and~(\ref{Gbound}), 
it follows that if $A \in \A^n$ and $x \in A^{*,n}$, 
then   $\dist(x,\bd A)
 \geq c  n^{7/8}$. That is, if we translate $A$ to make $x$ the origin, 
then the inradius of the translated set is at least $c n^{7/8}$.
  Hence, we can use Theorem~\ref{greentheoremB} to deduce~(\ref{mar29.eq2}).
\end{proof}

\subsection{An estimate for hitting the boundary} \label{hit-section}

Suppose that $A$ is any finite, connected subset of $\Z^2$, not 
necessarily simply connected, and let $V \subset \bd A$
be non-empty. Recall that $\tilde A \subset \C$ is the ``union of
 squares'' domain associated to $A$ as in Section~\ref{scs}. For
 every $y \in V$, consider the collection of edges containing $y$, 
namely $\edgeset{y} := \{(x,y) \in \bd_e A\}$, and set $\edgeset{V} := 
\bigcup_{y\in V} \edgeset{y}$. If $\ell_{x,y}$ is the perpendicular line 
segment of length one intersecting $(x,y)$ in the midpoint as in Section~\ref{scs}, 
then define 
$\tilde V := \bigcup\limits_{(x,y) \in \edgeset{V}} \ell_{x,y}$
to be the associated boundary arc in $\bd \tilde{A}$. Suppose that
 $\tau_A := \min\{j : S_j \in \bd A\}$, $T_A = T_{\tilde A} := 
\inf\{t : B_t \in \bd \tilde{A}\}$, and throughout this subsection, write
\begin{equation}\label{h-def}
h(x) = h_A(x,V) :=  \P^x\{S_{\tau_A} \in V \} = \sum_{y\in V} h_A(x,y),
\end{equation}
and
\begin{equation}\label{H-def}
H(x)  = H_{\tilde A}(x,\tilde{V}) :=  \P^x\{B_{T_A} \in \tilde V\} = \int_{\tilde V} 
H_{\tilde A}(x,y) \,|\d y|,
\end{equation}
where $h_A$ and $H_{\tilde A}$ are the discrete Poisson kernel and the Poisson kernel, 
respectively.

\begin{definition}
If $F:\Z^2 \to \R$, let $L$ denote the \defines{discrete Laplacian} defined by
$$LF(x) :=\E^x[F(S_1)-F(S_0)]
 =    \frac{1}{4}\sum_{|e|=1}(F(x+e)-F(x)).$$
Call a function $F$ \defines{discrete harmonic at $x$} if $LF(x)=0$.  
If $LF(x)=0$ for all $x\in A \subseteq \Z^2$, then $F$ 
is called \defines{discrete harmonic in $A$}.
\end{definition}

Let $\Delta$ denote the usual Laplacian in $\C$, and recall that $F:\C
\rightarrow \R$ is
 \defines{harmonic at $x$} if $\Delta F(x) = 0$. Note that $L$ is a natural 
discrete analogue of $\Delta$. 
 If $r>1$, $F:\C \to \R$, and $\Delta F(x)=0$ for 
all $x \in \C$ with $|x|<r$, then   Taylor's series and uniform 
bounds on the derivatives of harmonic functions~\cite{Bass} imply that
\begin{equation}\label{sept13.1}
|LF(0)| \le \|F\|_{\infty} O(r^{-3}).
\end{equation}
Note that $h$ defined by discrete~(\ref{h-def}) is  
  discrete harmonic in $A$,
and    $H$ defined by~(\ref{H-def}), is harmonic in $\tilde{A}$. 
Our goal in the remainder of this subsection is to prove the following 
proposition.

\begin{proposition}\label{aug30.1}
For every $\eps>0$, there exists a $\delta > 0$ such that if $A$ is a finite 
connected subset of $\Z^2$, $V \subset \bd A$, and $x \in A$ with $H(x) \ge 
\epsilon$, then $h(x) \ge \delta$.
\end{proposition}

We first note that for every $n < \infty$, there is a $\delta' = \delta'(n) > 0$ such
 that the proposition holds for all $A$ of cardinality 
at most $n$ and all $\eps > 0$. This is because $h$ and $H$ are 
strictly positive (since $V \neq \emptyset$) and the collection of connected 
subsets of $\Z^2$ containing the origin of cardinality at most $n$ is finite.  Hence, 
we can choose
\begin{equation}\label{deltaprime}
\delta'(n) = \min \P^x\{S_{\tau_A} = y\}
\end{equation}
where the minimum is over all finite connected $A$ of cardinality at most
$n$, all $x \in A$, and all $y \in \bd A$. We now extend this to all $A$ for 
$x$ near the boundary. 

\begin{lemma}\label{aug30.2}
For every $\eps > 0  ,M < \infty$, there exists a $\delta > 0$, 
such that if $A$ is a finite connected subset of $\Z^2$,  $V \subset \bd A$,
 and $x \in A$ with $H(x) \ge \eps$ and $\dist(x,\bd \tilde A) \leq M$, then 
$h(x) \ge \delta$.
\end{lemma}

\begin{proof} 
By the recurrence of planar Brownian motion, we can 
find an $N = N(M,\eps)$ such that  
$\P^x\{ \diam(B[0,T_A]) \le N \} \ge 1-(\epsilon/2)$
whenever
$\dist(x,\bd \tilde A) \le M$.
Hence, if $\P^x\{B_{T_A} \in \tilde V\} \ge \epsilon$ and $\dist(x,\bd \tilde A)
 \le M$, then
$$\P^x\{B_{T_A} \in \tilde V; \diam(B[0,T_A])  \le N \} \ge \epsilon/2,$$
and the lemma holds with $\delta = \delta' (3N)$, say, where $\delta'$ is
 defined as in~(\ref{deltaprime}) above.
\end{proof}

For every $M < \infty$ and finite $A \subset \Z^2$, let
$$\sigma_M = \sigma_{M,A} := \min\{j \ge 0: \dist(S_j,\bd A) \le M \}.$$
Since $A$ is finite and $h$ is a discrete harmonic function on $A$, it is necessarily 
bounded so that $h(S_{n \wedge \sigma_M})$ is a bounded martingale. It then 
follows from the optional sampling theorem that $h(x) = \E^x[h(S_{\sigma_M})]$
 for all $x \in A$ since $\sigma_M$ is an a.s.~finite stopping time. The next 
lemma gives a bound on the error in this equation if we replace $h$ with $H$.

\begin{lemma}\label{apr9.lem1}  
For every $\eps>0$, there exists an $M<\infty$ such that if $A$ is a finite 
connected subset of $\Z^2$, $V \subset \bd A$, and $x \in A$, 
then $|\, H(x) -  \E^x[H(S_{\sigma_M})] \, | \le \eps$.
\end{lemma}

\begin{proof} 
For any function $F$ on $A$ and any $x \in A$,
\begin{equation}\label{apr9.eq1}
F(x) = \E^x[F(S_{\sigma_M})] - \E^x[\sum_{j=0}^{\sigma_M - 1} LF(S_j) ].
\end{equation}
(Note that $F$ is bounded since $A$ is finite.) Applying~(\ref{apr9.eq1}) to $H$ gives
$$|\: H(x) -  \E^x[H(S_{\sigma_M})] \: |  \le \sum_{\dist(y,\bd A) \ge M}
 G_A(x,y) \: |LH(y)|.$$ 
Since $H$ is harmonic and bounded by 1,~(\ref{sept13.1}) implies $|LH(y)| 
\le c \: \dist(y,\bd A)^{-3}$. A routine estimate shows that there is a 
constant $c$ such that for all $A \in \A$, $x \in A$, and $r \ge 1$, if a
simple random walk is within distance $r$ of the boundary, then the probability
 that it will hit the boundary in the next $r^2$ steps is bounded below by
 $c/\log r$. Consequently, we have
$$\sum_{r \le \dist(y,\bd A) \le 2r} G_A(x,y) \le c \: r^2 \log r.$$
Combining these estimates gives
$$\sum_{r \le \dist(y,\bd A) \le 2r} G_A(x,y) \: |L H(y)| \le c r^{-1} \log r,$$
and hence
\begin{equation}\label{above}
\sum_{\dist(y,\bd A) \ge M} G_A(x,y) \: |L H(y)| \le  c\;  M^{-1} \; \log M.
\end{equation}
The proof is completed by choosing $M$ sufficiently large which will 
guarantee that the right side of~(\ref{above}) is smaller than $\eps$.
\end{proof}
  
\begin{proof}[Proof of Proposition~\ref{aug30.1}]
Fix $\epsilon>0$, and suppose $H(x) \geq \epsilon$.  By the Lemma~\ref{apr9.lem1}, 
we can find an $M = M(\epsilon)$ such that
$$\E^x[H(S_{\sigma_M})] = \sum_{\dist(y,\bd  A) \le M} J(x,y) \:H(y) \ge \epsilon/2,$$
where $J(x,y) =  J_A(x,y;M) := \P^x\{S_{\sigma_M} = y\}$.  Hence
$$\sum_{H(y) \ge \eps/4, \: \dist(y,\bd  A)\le M} J(x,y) \; 
\ge \; \sum_{H(y) \ge \eps/4, \: \dist(y,\bd A)\le M} J(x,y) \:H(y) \ge \eps/4.$$
By Lemma~\ref{aug30.2}, there is a $c = c(\epsilon,M)$ such that $h(y) \ge c$ if 
$H(y) \ge \eps/4$ and $\dist(y,\bd  A) \le M$. Hence,
\begin{equation*}
h(x) = \sum_{\dist(y,\bd A)\le M} J(x,y) \: h(y) \ge  c\: \eps/4. \qedhere
\end{equation*}
\end{proof}

\subsection{The main estimates}\label{boundary-sect}

Theorem~\ref{greentheoremB} with Corollary~\ref{sept9.prop2.nowcor} give good estimates if $f_A(x)$ and $f_A(y)$ 
are not close to $\bd \D$. While the result is true even for points 
near the boundary, it is not very useful because the error terms are much 
larger than the value of the Green's function. Indeed, if $A \in \A^n$ and 
$x\in \bd_i A$, then $g_A(x) = O(n^{-1/2})$ and $G_A(x) =O(n^{-1/2})$, but
 $O(n^{-1/2}) \ll O(n^{-7/24}\log n)$, the error term in Corollary~\ref{sept9.prop2.nowcor}.

 In this subsection we establish Proposition~\ref{sept23.prop1} which gives 
estimates for $x$ and $y$ close to the boundary provided that they are not
 too close to each other.  Theorem~\ref{rwestimate}
follows immediately. Recall that $A^{*,n} := \{ x \in A: g_A(x)
 \geq n^{-1/16} \}$ as in~(\ref{stardomain}).

The following estimates can be derived easily from~(\ref{sept9.1}). If $z=f_A(x)=(1-r)
 e^{i\theta}$, $z'=f_A(y) \in \Disk$ with $|z-z'|\geq r$, then
\begin{equation}\label{oct30.1}
g_A(x,y) = g_{\D}(z,z') = \frac{g_{\D}(z)\;(1-|z'|^2)}{|z'-e^{i \theta}|^2}\; 
\left[1+O\left(\frac{r}{|z-z'|}\right)\right].
\end{equation}
Similarly, if $z'=f_A(y)=(1-r') e^{i\theta'}$ with $r \geq r'$ and $|z-z'| \geq r$, 
\begin{equation}\label{oct30.2}
g_A(x,y) = g_{\D}(z,z') = \frac{g_{\D}(z)\; g_{\D}(z')}{1-\cos(\theta-\theta')}\;
\left[1+O\left(\frac{r}{|\theta-\theta'|}\right)\right]. 
\end{equation}
 
\begin{proposition}\label{sept23.prop1}
Suppose $A \in \A^n$. If $x \in A \setminus A^{*,n}$ and 
\begin{equation}\label{Jdomain}
J_{x,n} := \{z \in A: |f_A(z)-\exp\{i\theta_A(x)\} | \geq n^{-1/16} \log^2 n \},
\end{equation}
then for $y \in J_{x,n}$,
\begin{equation}\label{oct30.3}
G_A(x,y)=G_A(x)\; \frac{1-|f_A(y)|^2}{|f_A(y)-e^{i \theta_A(x)}|^2}\; 
\left[\;1+O\left(\frac{n^{-1/16} \log n}{|f_A(y)- e^{i\theta_A(x)}|}\right)\;
\right], \;\;\;\; y \in A^{*,n},
\end{equation}
\begin{equation}\label{oct30.4}
G_A(x,y)=\frac{\;(\pi/2)\; G_A(x)\; G_A(y)}{1-\cos(\theta_A(x)-\theta_A(y))}\; 
\left[\;1+O\left(\frac{n^{-1/16}\log n}{|\theta_A(y)-\theta_A(x)|}\right)\;\right], \;\;\;\;
 y\in A\setminus A^{*,n}.
\end{equation}
\end{proposition}

Thus, in view of the estimates~(\ref{oct30.1}) and~(\ref{oct30.2}) above, 
there is nothing surprising about the leading terms in~(\ref{oct30.3}) 
and~(\ref{oct30.4}). Proposition~\ref{sept23.prop1} essentially says that
 these relations are valid, at least in the dominant term, if we replace $g_A$ 
with $(\pi/2)\: G_A$.  Theorem~\ref{rwestimate} follows from this
proposition and~(\ref{dec15.1}).

The hardest part of the proof is a lemma that states that if the random walk 
starts at a point $x$ with $f_A(x)$ near $\bd \D$, then, given that the walk 
does not leave $A$, $f_A(S_j)$ moves a little towards the centre of the disk 
before its argument changes too much.

\begin{lemma}\label{sept10.lemma2}
For $A \in \A^n$, let
$\eta = \eta(A,n) := \min\{j \geq 0: S_j \in A^{*,n} \cup A^c\}$.
There exist constants $c$, $\alpha$ such that if $A \in \A^n$, $x\in A 
\setminus A^{*,n}$, and $r > 0$, then
$$\P^x\{\max_{0 \leq j \leq \eta-1}|f_A(S_j)-f_A(x)|\geq r \: n^{-1/16}\} 
\leq c \:e^{-\alpha r}\:\P^x\{S_\eta \in A^{*,n}\},$$
$$\P^x\{|f_A(S_\eta)-f_A(x)|\geq r \:n^{-1/16} \mid S_\eta \in A^{*,n} \}
 \leq c \: e^{-\alpha r}.$$
In particular, there is a $c_0$ such that if
$$\xi =\xi(A,n,c_0) := \min\{j \ge 0:S_j \not\in A
\mbox{ or }  |f_A(S_j)-f_A(x)|\geq c_0\; n^{-1/16} \; \log n \},$$
then 
\begin{equation*}
\P^x\{\xi < \eta\} \leq c_0 \; n^{-5}\; \P^x\{S_\eta \in A^{*,n}\}.
\end{equation*}
\end{lemma}
 
In order to prove this lemma, we will need to establish several ancillary results. Therefore, we devote Section~\ref{harnsec} which follows to the complete proof of this lemma, and then prove Proposition~\ref{sept23.prop1} in the separate Section~\ref{prop-proof}.

\subsubsection{Proof of Lemma~\ref{sept10.lemma2}}\label{harnsec}

If $A \in \A$ and $x \in A$, let $d_A(x)$ be the distance from $f_A(x)$ to the unit circle. Note that
$d_A(x) = \dist(f_A(x), \bd \D) = 1 - |f_A(x)| = 1-\exp\{-g_A(x)\}$ in view of~(\ref{oct5.1}).
As a first step in proving Lemma~\ref{sept10.lemma2}, we need the following.

\begin{lemma}\label{sept10.lemma3}
There exist constants $c$, $c'$, $c''$, $\eps$ such that if $A \in \A$, $x \in A$ with $d_A(x) \leq c$, and $\sigma = \sigma(x,A,c,c')$ is defined by 
$$\sigma := \min\{j \geq 0 : S_j \notin A, \; d_A(S_j) \geq (1+c) d_A(x), \text{ or } |\theta_A(S_j) - \theta_A(x)| \geq c' d_A(x)\},$$ 
then $\P^x\{S_\sigma \notin A\} \geq \eps$, $\P^x\{S_\sigma \in A ; \, d_A(S_\sigma) \geq (1+ c) d_A(x) \} \geq \eps$, and
\begin{equation}\label{apr9.eq2}
\P^x\{ |\theta_A(S_\sigma) - \theta_A(x)| \leq c'' d_A(x) \mid S_\sigma \in A \} = 1.
\end{equation}
\end{lemma}

\begin{remark}
(\ref{apr9.eq2}) is not completely obvious since the random walk takes discrete steps.
\end{remark}

\begin{proof}
We start by stating three inequalities
 whose verification we leave to the reader. 
These are simple estimates for conformal maps on domains 
that are squares or the unions of two squares. Recall that
 $\Square_x$ is the closed square of side one centred at $x$ whose
 sides are parallel to the coordinate axes. There exists a constant
 $c_2 \in (1, \infty)$ such that if $A 
\in \A$; $x$, $y$, $w \in A$; $x \neq y$; $|w-x| = 1$, 
 then $d_A(z) \leq c_2 d_A(x)$ for $z \in \mathcal{S}_x$; $|f_A(z)-f_A(x)|
 \leq c_2 \, |f_A(z')-f_A(x)|$ for $z$, $z' \in \mathcal{S}_y$; 
and $|f_A(w) - f_A(x)|\leq c_2 d_A(x).$
The first of these inequalities implies that if $z \in \mathcal{S}_y$ 
and $d_A(z) \geq 3 c_2 d_A(x)$, then $d_A(y) \geq 3 d_A(x)$.
Fix $A \in \A$, $x \in A$ with $d_A(x) \leq 1/(100c_2^2)$, and let 
$$J = J(x,A) := \big\{\, y \in A: \mathcal{S}_y \cap \{z \in \tilde A : |f_A(z) - f_A(x)| < 5 c_2 d_A(x) \} \neq \emptyset \,\big\}.$$
That is, $y \in J$ if there is a $z \in \Square_y$ with $|f_A(z)-f_A(x)|< 5 c_2 d_A(x)$. Note that $J$ is a connected subset of $A$ (although it is not clear whether it is \emph{simply} connected) and
$\tilde J \subseteq \{z \in \tilde A : |f_A(z)-f_A(x)| < 5 c_2^2 d_A(x) \}$.
In particular, $d_A(y) \leq 6 c_2^2 d_A(x)$
and $|\theta_A(y) - \theta_A(x)| < c'd_A(x)$ for all $y \in J$.  

There is a positive probability $\rho_1$ that 
a Brownian motion in $\D$ starting at $f_A(x)$ leaves
 $\D$ before leaving the disk of radius $2 d_A(x)$ about
 $f_A(x)$. By conformal invariance, this implies that with probability 
at least $\rho_1$, a Brownian motion starting at $x$ leaves $\tilde J$ 
at $\bd \tilde A$.  Hence by Proposition~\ref{aug30.1}, there is an $\eps_1$
 such that $\P^x\{ S_{\tau_J} \not\in A \} \geq \eps_1$. 

Similarly, there is a positive probability $\rho_2$ that 
a Brownian motion in the disk starting at $f_A(x)$ reaches 
$[1-6 c_2^2 d_A(x)] \, \D$ before  leaving $\D$ and before leaving the set
$\{z: \, |\theta_A(z) - \theta_A(x)| \leq d_A(x)\}$.
Note that $d_A(z)^2 \geq  |f_A(z) - f_A(x)|^2 - d_A(x)^2$ on this set. In particular,
 with probability at least $\rho_2$, a Brownian motion starting at $x$ leaves $\tilde J$ 
at a point $z$ with $d_A(z) \geq 4 c_2 d_A(x)$.  Such a point $z$ must
 be contained in an $\mathcal{S}_y$ with $d_A(y) \geq 4 d_A(x)$. Hence,
 again using Proposition~\ref{aug30.1}, 
there is a positive probability $\eps_2$ that
 a random walk starting at $x$ reaches a point $y \in J$ with $d_A(y)
 \geq 4 d_A(x)$ before leaving $J$. 
In the notation of Lemma~\ref{sept10.lemma3} choose $c:=1/(100c_2^2)$,
let $c'$ be the $c'$ mentioned above, and let $\eps := \min\{\eps_1, \eps_2\}$. 
 Then we have already shown that $\P^x\{S_\sigma \notin A \} \geq \eps$ and  
$\P^x\{S_\sigma \in A ; \, d_A(S_\sigma) \geq (1+ c) d_A(x) \} \geq \eps$. 
Also, if $y$, $w \in A$ with $|y-w| = 1$, $d_A(y) \leq d_A(x)$, 
and $|\theta_A(y) - \theta_A(x)|  \leq c' d_A(x)$, then 
$$|f_A(w)-f_A(x)| \leq |f_A(w)-f_A(y)|+|f_A(y)-f_A(x)| \leq c''' d_A(x),$$
which implies that $|\theta_A(w) - \theta_A(x)| \leq c'' d_A(x)$ for an
 appropriate $c''$. This gives the last assertion in Lemma~\ref{sept10.lemma3}
 and completes the proof.
\end{proof}

\begin{corollary}  \label{sept12.cor1}
There exist constants $c$, $c'$, $\alpha$ such 
that if $a \in (0,1/2)$,  $A \in \A$, and $x \in A$ with
  $d_A(x) < a$, then the probability that a random walk starting at $x$ 
reaches the set $\{y \in A: d_A(y) \geq a \}$ without leaving the
 set $\{y \in A:|\theta_A(y) - \theta_A(x)| \leq c' a\}$ is at least $c(d_A(x)/a)^\alpha$.
\end{corollary}

\begin{proof}
Let $Q(r,a,b)$ be the infimum over all $A \in \A$ and
 $x \in A$ with $d_A(x) \geq r$ of the probability that a random walk
 starting at $x$ reaches the set $\{y \in A: d_A(y) \geq a \}$  without
 leaving the set $\{y \in A:|\theta_A(y) - \theta_A(x)| \leq b \}$.  It
 follows from Lemma~\ref{sept10.lemma3} that there exist $q > 0$, $\rho < 1$,
 and a $c''$ such that $Q(\rho^k,\rho^j,b) \geq q \; Q(\rho^{k-1} ,\rho^j, b - c'' \rho^k)$.  By iterating this we get
$Q(\rho^k,\rho^j,2 \; c''\; (1-\rho)^{-1} \;   \rho^j) \geq q^{k-j}$.
This and obvious monotonicity properties give the result.
\end{proof}

\begin{remark}
 A similar proof gives an upper bound of $c_1( d_A(x)/a)^{\alpha_1}$, 
but this result is not needed.
\end{remark}

For any $a \in (0,1)$ and any $\theta_1 < \theta_2$, 
let $\eta(a,\theta_1,\theta_2)$ be the first
 time $t \geq 0$ that a random walk  leaves the set $\{y \in A: d_A(y) \leq
 a, \theta_1 \leq \theta_A(y) \leq \theta_2\}$. Let 
$$q(x,a,\theta_1,\theta_2) := \P^x \{d_A(S_{\eta(a,\theta_1,\theta_2)}) > 
a \mid S_{\eta(a,\theta_1,\theta_2)} \in A\},$$
and note that if $\theta_1 \leq \theta_1' \leq \theta_2' \leq 
\theta_2$, then $q(x,a,\theta_1',\theta_2')\leq  q(x,a,\theta_1,\theta_2)$.

\begin{proposition}\label{oct12.prop1}
There exist constants $c$, $c_1$ such that if $a \in (0,1/2)$,
 $A \in \A$, and $x \in A$, then   
$q(x,a,\theta_A(x) - c_1a, \theta_A(x) + c_1a) \geq c$.
\end{proposition}

\begin{proof}
For every $r>0$ and $m \in \N$ let 
$h(m,r) := \inf q(x,a,\theta_A(x) - ra, \theta_A(x) +ra)$ 
where the infimum is taken over all $a \in (0,1/2)$, $A \in \A$, and 
 all  $x \in A$ with  $d_A(x) \geq 2^{-m} a$. The proposition is
 equivalent to saying that there is a $c_1$ such that
$$\inf_m h(m,c_1) > 0.$$
It follows from Corollary~\ref{sept12.cor1} that there is a $c'$ 
such that $h(m,c') >0$ for each $m$; more specifically, 
there exist $c_2$, $\beta$ such that $h(m,r) \geq c_2 e^{-\beta m }$ for  $r \geq c'$.

Suppose $x \in A$ with $2^{-(m+1)}a \leq d_A(x) < 2^{-m}a$. Start a 
random walk at $x$ and stopped at $t^*$, defined to be the first time $t$ 
when one of the following is satisfied:   $S_t \notin A$, $d_A(S_t) \geq 2^{-m}a$, 
or $|\theta_A(S_t)-\theta_A(x) | \geq m^2 2^{-m} a$.  By 
iterating Lemma~\ref{sept10.lemma3}, we see that the probability 
that the last of these three possibilities occurs is bounded above 
by $c'' e^{-\beta' m^2 }$. Choose $M$ sufficiently large such that 
for $m \geq M$, the last term is less than $c_2 e^{-2 \beta (m+1) }$, and 
such that
\begin{equation}\label{apr9.eq3}
\P^x\{ | \theta_A(S_{t^*}) - \theta_A(x) | \leq  m^3  2^{-m} a \mid S_{t^*} \in A \}=1.
\end{equation}
Note that~(\ref{apr9.eq2}) shows~(\ref{apr9.eq3}) holds for all 
sufficiently large $m$; for such $m$, if $r \geq c'$, then
$h(m+1,r) \geq [1 - e^{-\beta m }]\;  h(m,r - m^3  2^{-m})$.
In particular, if 
$$r > c' + \sum_{m=1}^\infty m^3 2^{-m},$$
then
\begin{equation*}
h(m,r) \geq [ \prod_{j=1}^\infty   (1 - e^{-\beta j })] \;h(M,c') > 0. \qedhere
\end{equation*}
\end{proof}

\begin{corollary}\label{apr8.cor1}
There exist $c$, $\beta$ such that if $a \in (0,1/2)$, $r > 0$, $A \in \A$, 
and $x \in A$, then
$q(x,a,\theta_A(x) - ra, \theta_A(x) + ra) \geq 1 - c e^{-\beta r}$.
\end{corollary}

\subsubsection{Proof of Proposition~\ref{sept23.prop1}}\label{prop-proof}

Suppose $A \in \A^n$,  $(z,y) \in \bd_e A^{*,n}$, and
let $\eta$ be as in Lemma~\ref{sept10.lemma2}. Since $g_A$ is
 harmonic in the disk of radius $O(n^{7/8})$ about $z$, standard estimates 
for positive harmonic functions give 
$$|g_A(z) - g_A(y)| \leq O(n^{-7/8})\: g_A(z) \leq o(n^{-7/8}).$$
Since $g_A(y) < n^{-1/16} \leq g_A(z)$, we conclude
$g_A(z) =  n^{-1/16} + o(n^{-7/8})$,
and similarly for $g_A(y)$.  Hence, by Theorem~\ref{greentheoremB},
$$G_A(z) = (2/\pi)\; n^{-1/16} + O(n^{-1/3} \log n)
         = (2/\pi)\; n^{-1/16}\;[1 + O(n^{-13/48} \log n)],$$
and similarly for $G_A(y)$. Therefore, for any $x \in A \setminus A^{*,n}$, 
\begin{align}\label{sept13.3}   
G_A(x) &= \P^x\{S_\eta \in A^{*,n}\} \; \E^x[G_A(S_\eta)
 \mid S_\eta \in A^{*,n}]\nonumber\\
       &=(2/\pi)\; n^{-1/16}\; \P^x\{S_\eta \in A^{*,n}\} \; [1 + O(n^{-13/48} \log n)]. 
\end{align}
In a similar fashion, note that if $x$, $y \in A^{*,n}$, then $g_A(x,y) \geq c n^{-1/8}$,
 and hence by Corollary~\ref{sept9.prop2.nowcor}, if $|x-y|\geq n^{1/2}$, then
$G_A(x,y) = (2/\pi)\; g_A(x,y) \; [1 + O(n^{-1/6} \log n)]$.

If $A \in \A^n$, $x \in A \setminus A^{*,n}$, and $y \in A^{*,n}$, then
 the strong Markov property gives
$$G_A(x,y) = \sum_{z \in A^{*,n}} \P^x\{S_\eta = z\} G_A(z,y).$$
If $x \in A \setminus A^{*,n}$, let 
$R(x) = R(x,n,A) := \{z \in A: |f_A(z) - f_A (x)| \leq c_0 n^{-1/16} \log n \}$,
where $c_0$ is the constant in Lemma~\ref{sept10.lemma2}.  From that lemma we see that   
$$\sum_{z \in A^{*,n} \cap R(x)} \P^x\{S_\eta = z\} = [1 - O(n^{-5})] \;
\sum_{z \in A^{*,n}} \P^x\{S_\eta = z \}.$$
But $c\; n^{-1/8} \leq G_A(z,y) \leq c' \log n$ for $z$, $y \in A^{*,n}$; hence,
$$G_A(x,y) = [1 + o(n^{-4})] \; \sum_{z \in A^{*,n} \cap R(x)}
 \P^x\{ S_\eta = z\}\;G_A(z,y).$$
If $|f_A(y) - f_A(x)| \geq n^{-1/16} \; \log^2 n$, and $z \in A^{*,n} \cap R(x)$,
 then from~(\ref{oct30.1}),
$$g_A(z,y) = \frac{n^{-1/16} \;(1-|f_A(y)|^2)}{|f_A(y)-e^{i\theta_A(x)}|^2}\;
[\; 1 + O(\frac{n^{-1/16} \; \log n}{|f_A(y)-e^{i\theta_A(x)}|}) \;].$$
Hence, using Corollary~\ref{sept9.prop2.nowcor},
$$G_A(x,y)=\P^x\{S_\eta \in A^{*,n}\}\; \frac{(2/\pi)\; n^{-1/16}\; (1-|f_A(y)|^2)}
{|f_A(y)-e^{i\theta_A(x)}|^2}\; [\; 1 + O(\frac{n^{-1/16} \; \log n}
{|f_A(y)-e^{i\theta_A(x)}|})\; ].$$
Combining this with~(\ref{sept13.3}) gives~(\ref{oct30.3}). If $y \in \bd_i A^{*,n}$, 
then we can write
\begin{equation}\label{hallow.2}
G_A(x,y) = G_A(x)\frac{n^{-1/16}}{1-\cos(\theta_A(x)-\theta_A(y))}\; 
[\;1 + O(\frac{n^{-1/16} \log n}{|y - e^{i\theta_A(x)}|})\; ]. 
\end{equation}

Now suppose $y \in J_{x,n} \setminus A^{*,n}$. Then we can write
$$G_A(x,y) = G_{A \setminus A^{*,n}}(x,y)+ 
\sum_{z \in A^{*,n}} \P^x\{S_\eta = z \} G_A(z,y),$$
and using~(\ref{oct30.3}) on $G_A(z,y)$ gives 
$G_A(x,y) \geq c  n^{-1/8}  \P^x\{S_\eta \in A^{*,n}\}  \P^y \{S_\eta \in A^{*,n}\}$.  
However, provided $R(x) \cap R(y) =\emptyset$, which is 
true for $n$ sufficiently large, Lemma~\ref{sept10.lemma2} shows that
$G_{A \setminus A^{*,n}}(x,y)\leq c n^{-10} \P^x\{S_\eta \in A^{*,n}\}
 \P^y \{S_\eta \in A^{*,n}\}$.
Therefore,
$$G_A(x,y) = [1 + o(n^{-9})] \; \sum_{z \in A^{*,n}} \P^x\{ S_\eta = z \} G_A(z,y).$$
and we can use~(\ref{hallow.2}) to deduce~(\ref{oct30.4}).

\section*{Acknowledgements}
Much of this research was done by the first author
 in his Ph.D.~dissertation~\cite{KozThesis} under the
 supervision of the second author.  The first author expresses his
 gratitude to the second author for his continued guidance and support.


\bibliographystyle{plain}


\vfill

\begin{center}
\begin{tabular}{lll}
Michael J.~Kozdron                      &\hspace{1in}   &Gregory F.~Lawler\\
\texttt{kozdron@math.uregina.ca}        &   &\texttt{lawler@math.cornell.edu}\\
Department of Mathematics \& Statistics &   &Department of Mathematics\\
College West 307.14                     &   &Malott Hall\\
University of Regina                    &   &Cornell University\\
Regina, SK S4S 0A2                      &   &Ithaca, NY 14853-4201
\end{tabular}
\end{center}

\end{document}